\documentclass[a4paper]{article}

\usepackage{paperdraft}
\usepackage[draft]{changes}
\definechangesauthor[name=Jose, color=blue]{JOS}
\title{A Scalable Algorithm for Shape Optimization with Geometric Constraints in Banach Spaces}
\usepackage{booktabs}
\usepackage{subcaption}
\usepackage{cleveref}
\author[2]{P. Marvin M\"uller}
\author[1]{J. Pinz\'on}
\author[2]{T. Rung}
\author[1]{M. Siebenborn}
\affil[1]{Hamburg University}
\affil[2]{Hamburg University of Technology}

\date{\today}
\usepackage{biblatex}
\begin{document}

\maketitle

\begin{abstract}
    This work develops an algorithm for PDE-constrained shape optimization based on Lipschitz transformations.
    Building on previous work in this field, the $p$-Laplace operator is utilized to approximate a descent method for Lipschitz shapes.
    In particular, it is shown how geometric constraints are algorithmically incorporated avoiding penalty terms by assigning them to the subproblem of finding a suitable descent direction.
    A special focus is placed on the scalability of the proposed methods for large scale parallel computers via the application of multigrid solvers.
    The preservation of mesh quality under large deformations, where shape singularities have to be smoothed or generated within the optimization process, is also discussed.
    It is shown that the interaction of hierarchically refined grids and shape optimization can be realized by the choice of appropriate descent directions.
    The performance of the proposed methods is demonstrated for energy dissipation minimization in fluid dynamics applications.
\end{abstract}

\par
{\noindent\small \textbf{Keywords}: Shape optimization, Lipschitz transformations, p-Laplace, geometric multigrid, parallel computing}
\par

\section{Introduction}
\label{sec:introduction}
In this paper we present a numerical scheme for the efficient treatment of geometrical constraints in shape optimization within the context of fluid dynamics applications.
Several advantages over other well-known techniques are described, particularly how the need for penalty terms is relinquished in favor of a more robust approach. 
Additionally, we argument how the presented algorithm is well-suited for geometrical constraints of integral form, which are preserved up to a numerical tolerance during the optimization process.

Constraints on the volume and barycenter are often required in fluid dynamics. 
This is particularly true for the minimal drag problem of a free floating body, see~\cite{schulz2016, onyshkevych2020, mueller2021, pinzon2021}
Another example is in~\cite{allaire2004structural}, where a volume and perimeter constraints are considered for a structural optimization problem.
A constraint for minimum and maximum thickness is formulated in~\cite{allaire2020, geiser2021aggregated} in order to meet requirements stemming from the manufacturing process.
 
In order to preserve these constraints, we include them in the process of finding descent directions in Banach spaces.
This is especially challenging, since the geometrical constraints are of a different type than that of the state equation, i.e. the governing partial differential equation (PDE).
Meaning, that the flow field is characterized by the stationary, incompressible Navier-Stokes equations, which lead to a PDE constraint optimization problem over infinite dimensional Banach spaces.
Whereas, the geometrical constraints are given by a finite number of integral type constraints, independently of a finite element model (see \Cref{sec:model}).

We focus on a well-established benchmark problem in the field of shape optimization constrained by Navier-Stokes equations, where the shape of an obstacle located within a flow channel is to be optimized with respect to the drag generated over its surface.
In general, the optimization problem can be summarized as determining a geometry $\Omega \subset \R^d$ that minimizes a shape functional $J$. 
The functional represents a physical quantity, e.g. drag of an obstacle or the energy dissipation associated with the flow around an obstacle such that a fixed number of geometric constraints $g(\Omega) = 0 \in \R^m$ is also fulfilled.
There are two major challenges in this problem, one concerning admissibility and optimality, and another related to the regularity of the obstacle's shape. 
In addition, the aforementioned geometrical constraints have to preserved in order to avoid trivial and non-feasible solutions.

Many popular approaches rely on strategies to simultaneously update state variables and Lagrange multipliers of the constraints, for an overview see~\cite{arora1991}.
In other words, optimality and admissibility are established simultaneously within one iteration.
However, the optimization problems addressed here tend to be non-meaningful or even unsolvable, provided that the constraints are not precisely fulfilled.
For instance, if the barycenter of the obstacle is not fixed, then the object would leave the domain. 
In the same way the optimization procedure usually yields a trivial solution if the volume is not preserved, because then the obstacle would be contracted to a single point.
The second issue is the regularity of initial and optimized shapes.
On one hand, it might be essential that singularities can be represented in the optimal shapes, e.g. kinks and sharp edges.
On the other, it is necessary for the scalability of the algorithm  to apply multigrid methods as a grid-independent preconditioner for a Krylov subspace solver.
In addition, the discretization of the domain with a coarse grid must be able to adequately represent both the surface of the obstacle and far field boundaries on the base level.  
Therefore during the optimization process the descent directions and shape updates have to feature non-smooth characteristics.

A common approach in shape optimization is to map a reference domain $\Omega \subset \R^d$ with $d = 2$ or $d = 3$, to a perturbed domain $(\id + u)(\Omega) := \{ x + u(x) \in \R^d : x \in \Omega \}$ with $u \in W^{1,\infty}(\R^d, \R^d)$ such that $\id + u$ is  a Lipschitz homeomorphism, cf. \cite{sokolowski1992, zolesio2011, allaire2020}. 
Note that this does not require any parametrization of the geometry, e.g. like in a CAD description with NURBS-surfaces. 
For the shape deformation we follow \cite{deckelnick2021, ishii2005limits} and consider the steepest descent direction in $W^{1,\infty}$-topology with a $p$-Laplace relaxation and the deformation vector field $u$ then is the solution to a minimization problem.
In contrast here we consider a constraint optimization problem in order to take the geometric constraints into account.

The well-known fluid dynamic example for a minimal drag problem in \cite{pironneau1973optimum} considers the volume constraint, which is one dimensional.  
As the shape update is performed solely in the surface normal direction, the corresponding Lagrange multiplier is given by the mean value of the deformation.
In \cite{mohammadi2004}, constraints as maximum thickness and volume have been taken into account via penalization of the cost function. 
Even though this approach allows also for more general shape deformations the shape update has to be rather small in order to keep the procedure numerically stable. 
Besides this one, two other approaches have successfully been applied to shape optimization problems of this kind. 
Firstly, an augmented Lagrange method can be used to determine the Lagrangian multipliers associated with the geometric constraints, \cite{allaire2004structural,schulz2016,allaire2020,mueller2021}. 
However, this approach has difficulties that are challenging to overcome, e.g. several parameters are problem dependent and have to be assigned to appropriate initial values.  
Furthermore, the constraints first have to be violated in order to determine the desired multipliers and the whole shape optimization problem has to be solved repetitively until the multipliers converge.
This can lead to unfeasible shapes throughout the optimization procedure, as was previously mentioned.
Secondly, the method of mappings \cite{onyshkevych2020,haubner2021} enables for the fulfillment of the geometric constraints up to machine accuracy. 
Therefore, the shape optimization problem only has to be solved once and within each iteration only feasible geometries are computed. 
Here the state as well as the adjoint variables are determined on the transformed domain by applying the perturbation of identity.
Thus, the descent vector field couples to all the constraints. 
Nevertheless, it is necessary to solve the fully coupled optimality system as a whole which is challenging, not only from an implementation point of view, but because it is computationally expensive.

In contrast, in this approach we use a second order method only to determine the steepest descent direction in a first order shape optimization. 
This reduces the dimensions of the linear systems to be solved.
However, our otpimization scheme requires that the geometric constraints only depend on the descent vector field $u$ and not on the physical state variables, e.g. velocity or pressure.
Therefore, the optimality system can be solved sequentially starting with the state, then the adjoints to the state, and finally the descent direction. 
This gives us the ability to handle problems with very large degrees of freedom (DoFs), while fulfilling the geometric constraints up to an arbitrary tolerance. 
It may be mentioned that the approach presented has analogies to optimization on manifolds and the investigations carried out in \cite{schiela2019composite,schiela2021sqp}, although it is based on Hilbert space settings.
We also want to mention the Uzawa iteration \cite{uzawa1958iterative} for solving a saddle point problem which occurs in the presented approach.

The remainder of this paper has the following structure: In \cref{sec:model}, the physical problem is introduced and the basics of shape optimization in Banach spaces are recalled.
\Cref{sec:algorithm} proposes an algorithm which determines admissible shape deformation descent directions.
In \cref{sec:numerics}, we demonstrate a scalable multigrid implementation for a fluid dynamics benchmark problem, while the performance of the method is investigated in \cref{sec:scalability}.
In \cref{sec:conclusion} the presented algorithm and numerical experiments are recalled and summarized.

Regarding the notation in the upcoming equations $D(\cdot)$ denotes the Jacobian, for the spatial Euclidean gradient operator we use $\nabla (\cdot)$ and the directional derivatives with respect to a specific variable are indicated via, for instance $\frac{\partial}{\partial u} (\cdot) \test u$ in direction $\test u$.
The shape derivative of the functional $J(\Omega)$ in direction $u$ is denoted by $J'(\Omega) u$ as defined in \cref{eq:ShapeDerivative}.

\section{Model Equations}
\label{sec:model}

\begin{figure}[htp]
    \centering
    \def\svgwidth{0.8\textwidth}
\begingroup%
  \makeatletter%
  \providecommand\color[2][]{%
    \errmessage{(Inkscape) Color is used for the text in Inkscape, but the package 'color.sty' is not loaded}%
    \renewcommand\color[2][]{}%
  }%
  \providecommand\transparent[1]{%
    \errmessage{(Inkscape) Transparency is used (non-zero) for the text in Inkscape, but the package 'transparent.sty' is not loaded}%
    \renewcommand\transparent[1]{}%
  }%
  \providecommand\rotatebox[2]{#2}%
  \newcommand*\fsize{\dimexpr\f@size pt\relax}%
  \newcommand*\lineheight[1]{\fontsize{\fsize}{#1\fsize}\selectfont}%
  \ifx\svgwidth\undefined%
    \setlength{\unitlength}{281.86157947bp}%
    \ifx\svgscale\undefined%
      \relax%
    \else%
      \setlength{\unitlength}{\unitlength * \real{\svgscale}}%
    \fi%
  \else%
    \setlength{\unitlength}{\svgwidth}%
  \fi%
  \global\let\svgwidth\undefined%
  \global\let\svgscale\undefined%
  \makeatother%
  \begin{picture}(1,0.51746909)%
    \lineheight{1}%
    \setlength\tabcolsep{0pt}%
    \put(0,0){\includegraphics[width=\unitlength,page=1]{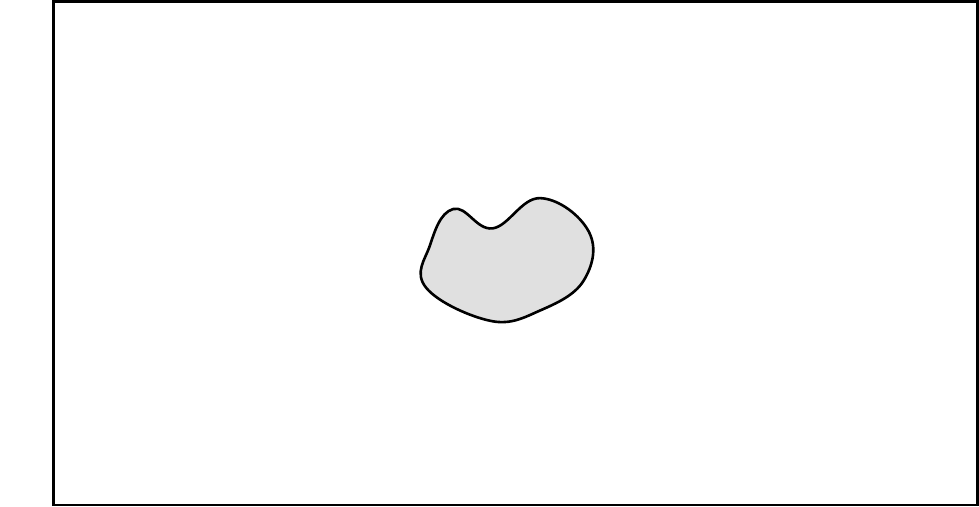}}%
    \put(-0.00294498,0.2670112){\makebox(0,0)[lt]{\lineheight{0}\smash{\begin{tabular}[t]{l}$\Gin$\end{tabular}}}}%
    \put(0.92832675,0.25149269){\makebox(0,0)[lt]{\lineheight{0}\smash{\begin{tabular}[t]{l}$\Gout$\end{tabular}}}}%
    \put(0.45113145,0.47723075){\makebox(0,0)[lt]{\lineheight{0}\smash{\begin{tabular}[t]{l}$\Gwall$\end{tabular}}}}%
    \put(0.44880138,0.02229574){\makebox(0,0)[lt]{\lineheight{0}\smash{\begin{tabular}[t]{l}$\Gwall$\end{tabular}}}}%
    \put(0.59605691,0.30603575){\makebox(0,0)[lt]{\lineheight{0}\smash{\begin{tabular}[t]{l}$\Gobs$\end{tabular}}}}%
    \put(0.71168408,0.13623505){\makebox(0,0)[lt]{\lineheight{0}\smash{\begin{tabular}[t]{l}$\Omega$\end{tabular}}}}%
    \put(0,0){\includegraphics[width=\unitlength,page=2]{domain.pdf}}%
    \put(0.46781223,0.22615455){\makebox(0,0)[lt]{\lineheight{0}\smash{\begin{tabular}[t]{l}$\Oobs$\end{tabular}}}}%
  \end{picture}%
\endgroup%

    \caption{Schematic view on a flow tunnel-like domain  $\Omega$ with the obstacle $\Oobs$ encircled by its surface $\Gobs$, wall boundaries $\Gwall$, inflow $\Gin$, and outflow $\Gout$. The height of the flow tunnel is defined by $\gamma>0$.}
    \label{fig:domain}
\end{figure}
In the present work, we propose an optimization methodology for PDE constraint shape optimization problems of the abstract form
\begin{align}
    \min\limits_{\Omega \in \mathcal{S}} &\quad j(\Omega,y) \label{eq:generalShapeOptProp1}\\
    \text{s.t.} &\quad e(\Omega,y) = 0 \label{eq:generalShapeOptProp2}\\
    &\quad g(\Omega) = 0 \label{eq:generalShapeOptProp3}
\end{align}
where $e$ denotes the PDE constraint on a state variable $y$.
The mapping $g$ refers to some finite dimensional geometric constraints on the bounded Lipschitz domain $\Omega$ with boundary $\Gamma = \Gobs \cup \Gin \cup \Gout \cup \Gwall$, where $\Gobs$ is to be optimized.
Furthermore, $\mathcal{S}$ denotes an abstract set of admissible shapes, as explained for instance in~\cite{sokolowski1992,zolesio2011}.

For some $m\in\mathbb{N}$ the geometric constraint is thus given as $g: \mathcal{S} \to \R^m$.
We assume the existence of the mapping $\Omega \mapsto y(\Omega)$. 
Thereby, we obtain the reduced cost functional $J(\Omega) := j(\Omega,y(\Omega))$.
In order to obtain sensitivities of the objective $J$ we follow \cite{sokolowski1992, zolesio2011,allaire2004structural,brandenberg2011}.
For this purpose, the domain $\Omega$ is parameterized in the sense of the perturbation of identity with the displacement field $u : \R^d \to \R^d$. 
For a sufficiently small $u$ we thus obtain deformed configurations
\begin{equation}
    \tilde{\Omega} := \left\{ x + u : x \in \Omega \right\}
    \label{eq:DefDeformedDomain}
\end{equation}
of the reference shape $\Omega$. For the sake of readability we abbreviate the perturbation of identity as
\begin{equation}
    F:\Omega \to \tilde{\Omega}; \quad F := \id + u. \label{eq:FDefinition}
\end{equation}
Moreover, the previously mentioned abstract set of admissible shapes $S$ can be further specified as 
\begin{equation}
	S:=\left\{F(\Omega):F =id+u, u\in W^{1,\infty}(\R^{d},\R^d)\right\}
\end{equation}

As parts of the boundary of $\Omega$ shall remain fixed, the displacement field $u$ is chosen in such a way that it vanishes on all boundaries of $\Omega$ which are not to be optimized. The directional shape derivative of $J$ evaluated in $\Omega$ in the direction $u$ is then given in a natural way by
\begin{equation}
    J(\tilde\Omega) = J(\Omega)%
    + J^\prime(\Omega) \, u
    + o(\| u \|) \quad \text{where}\quad \frac{o(u)}{\| u \|} \overset{\| u \| \to 0}{\longrightarrow} 0.
    \label{eq:ShapeDerivative}
\end{equation}
Then we can interpret the shape optimization problem \cref{eq:generalShapeOptProp1,eq:generalShapeOptProp2,eq:generalShapeOptProp3} locally, as a problem in $U_{ad} \subseteq W^{1,\infty}(\R^d, \R^d)$. 
In the present work, $U_\mathrm{ad}$ is the admissible set of displacements $u$ defining the transformation $F$,  which inherently fulfills the geometric constraints \cref{eq:generalShapeOptProp3} in the sense that $g(F(\Omega))=0$. 
Note that $U_\mathrm{ad} \neq \emptyset$ as $u = 0$ is an admissible transformation if the initial geometry $\Omega$ fulfills the geometrical constraints, i.e.\ $g(\Omega) = 0$.
The crucial aspect of the present method, is to separate the geometric constraints \cref{eq:generalShapeOptProp3} from the remaining PDE-constrained shape optimization problem \cref{eq:generalShapeOptProp1,eq:generalShapeOptProp2}, and move it to the admissible set $U_\mathrm{ad}$ of descent directions.
In contrast to other popular approaches, where admissibility is only guaranteed in the optimal configuration, we ensure that \cref{eq:generalShapeOptProp3} is fulfilled in each optimization step. This can be done, because the geometric constraints do not depend on the state $y$ and rely on the properties of the shape only.

For the computation of the shape derivative $J^\prime(\Omega)\, u$, we utilize formally the method of C\'ea, see for instance \cite[section~4.6]{allaire2020}.
In the following, we consider the particular problem of minimizing energy dissipation of the fluid flow mainly caused by an obstacle in a laminar, stationary flow, where the function 
\begin{equation}
    J(\Omega) := j(\Omega, v) = \frac{\nu}{2} \int_\Omega D \vel : D \vel \, dx
    \label{eq:EnergyDispObjectiveFunctional}
\end{equation}
is to be minimized.
Here and in the following, we denote velocity $\vel$, the density-specific pressure $\press$, viscosity $\nu$, and an inflow velocity $\vel_\infty$.
As the PDE constraint $e$, with state variable $y = (\vel,\press)$, we consider the stationary, incompressible Navier-Stokes equations
\begin{equation}
\begin{aligned}
    - \nu \Delta \vel + (\vel\cdot \nabla)\vel + \nabla \press = 0 & \text{ in } \Omega\\
    \Div{\vel} = 0  & \text{ in } \Omega\\
    \vel = 0  & \text{ on } \Gobs\cup\Gwall\\
    \vel = \vel_\infty  & \text{ on } \Gin\\
    D\vel \,\cdot n = \press n & \text{ on } \Gout.
    \end{aligned}
\label{eq:NavStokes}
\end{equation}
Here, we consider the adjoint approach for determining the directional shape derivative $J^\prime(\Omega)\, u$. For details on the adjoint Navier-Stokes equations, see e.g.~\cite{hinze2001,ulbrich2003,onyshkevych2020,pinzon2021}. For details on the shape derivative $J^\prime(\Omega)u$ of the objective function in \cref{eq:EnergyDispObjectiveFunctional}, see e.g. \cite{mohammadi2010,onyshkevych2020}.
In order to approximate the steepest descent direction $u$ in
\begin{equation}
    V^\infty_0 := \left\{ u \in W^{1,\infty}(\Omega \cup \Omega_\mathrm{obs},\R^d) : \| u \|_{W^{1,\infty}(\R^d, \R^d)} < 1, u = 0 \text{ a.e. } \Gin\cup\Gout\cup\Gwall \right\},
    \label{eq:V-0-infty}
\end{equation}
following \cite[Proposition 4.1]{allaire2020}, 
we introduce a $p$-Laplace relaxation with $p > 2$ with the corresponding minimization problem
\begin{equation}
\begin{aligned}
    \min_{u \in V^\infty_0} &\quad J^\prime(\Omega) \, u
    \\
    \text{s.t.} &\quad g(F(\Omega)) = 0,\\
                &\quad F = \id + u, 
\end{aligned} \label{eq:DefSteepestDescentDirection}
\end{equation}
inspired by \cite{ishii2005limits} and \cite{deckelnick2021}. 
This makes $F$ a Lipschitz transformation, not only on the obstacle's surface $\Gamma_{obs}$, but also on the entire domain $\Omega$.
Thus, in a discretization the elements undergo deformations that preserve mesh quality, as illustrated later in \cref{sec:numerics}.
Therefore, let
\begin{equation}
    V^p_0 = \left\{ u \in W^{1,p}(\Omega, \R^d): \|D u\|_{L^p(\Omega,\R^d)} \leq 1, u = 0 \text{ a.e. on } \Gin\cup\Gout\cup\Gwall \right\}
    \label{eq:V-0-p}
\end{equation}
and consider
\begin{equation}
\begin{aligned}
    \min_{u \in V^p_0} &\quad \frac{1}{p} \int_{\Omega} (Du : Du)^{p/2} \, dx + J'(\Omega)\, u\\
    \text{s.t.} &\quad g(F(\Omega)) = 0, \\
                &\quad F = \id + u
\end{aligned}
\label{eq:DefRelaxedSteepestDescentDirection}
\end{equation}
where it is assumed that $g(F(\Omega)): V^p_0 \to \R^m$, $u\mapsto g((\id +u)(\Omega))$, $m \geq 1$.
Notice that this is consistent with \cref{eq:generalShapeOptProp1,eq:generalShapeOptProp2,eq:generalShapeOptProp3} for $g$ over a fixed $\Omega$ and a variable displacement field $u$. 
In order to address the inconsistency between \cref{eq:V-0-infty} and \cref{eq:V-0-p}, we additionally assume that all displacements $u \in V_0^p$ are sufficiently small.  
Thus, the admissible set $\mathcal{S}$ is locally parameterized by $V^p_0$-deformations of $\Omega$. 
In the present work, $m=d+1$ refers to the barycenter and volume constraints
\begin{align}
    \int_{\Omega} (x + u) \det(DF) \; dx &= 0 \label{eq:BarycenterConstraint},\\
    \int_{\Omega} \det(DF) - 1 \; d x &= 0 \label{eq:VolumeConstraint},
\end{align}
and without loss of generality, we assume that the barycenter of the initial domain $\Omega$ is located at the origin $0\in \R^d$ of the domain, cf.~\cite{pinzon2021}. 
For computational reasons, the mathematical domain is restricted to the wetted domain, i.e. $\Omega$ without the shape $\Oobs$ itself. 
Furthermore, due to the constant volume constraint \cref{eq:VolumeConstraint}, we can omit in \cref{eq:BarycenterConstraint} the division by the reference volume and deformed domain $\Omega$ and $F(\Omega)$, respectively.
For the derivation of the optimality conditions of the steepest descent problem \cref{eq:DefRelaxedSteepestDescentDirection}, we define the Lagrangian function
\begin{equation}
\begin{aligned}
    L(u,\lambda) = &
    \frac{1}{p} \int_{\Omega} (Du : Du)^{p/2} \, dx + J'(\Omega)\, u \\
    &+ \sum_{i = 1}^d \lambda_i \int_{\Omega} (x_i + u_i) \det(DF) \; d x
    + \lambda_{d+1} \int_{\Omega} \det(DF) - 1 \; d x
\end{aligned}
\label{eq:Lagrangian}
\end{equation}
with $\lambda = (\lambda_1, \ldots, \lambda_d, \lambda_{d+1})^T$, where $\lambda_1, \ldots, \lambda_d$ are associated with the barycenter \cref{eq:BarycenterConstraint} and $\lambda_{d+1}$ with the volume constraint \cref{eq:VolumeConstraint}.
In the following we want to recall some rules of differentiation. 
Therefore, let $\test{u}, \stest{u}: \Omega \to \R^d$, $B: \Omega \to \R^{d\times d}$ and
\begin{equation}
    DF =\left(\frac{\partial}{\partial x_j} F_i\right)_{1\leq i,j \leq d} = I + Du
\end{equation}
the Jacobian of $F$. 
We specify the following useful formulae by applying the product and chain rule:
\begin{equation}
\begin{aligned}
    \frac{\partial}{\partial u} DF \; \test{u} &= D\test{u},\\
    \frac{\partial}{\partial u} \det(DF) \; \test{u} &= \tr(\invDF D\test{u})\det(DF),\\
    \frac{\partial}{\partial u} \left( \tr(DF \, B) \right) \; \test{u} &= B^T : \left( \frac{d}{du} DF \, \test{u} \right) = B^T : D\test{u},\\
    \frac{\partial}{\partial u} \left(\invDF\right) \; \test{u} &= - \invDF D\test{u} \invDF,\\
    \frac{\partial}{\partial u} \left( \tr( \invDF D\test{u} ) \right) \; \stest{u} &= 
    - D\test{u}^T : \invDF D\stest{u} \invDF\\ 
    &= \tr(-D\test{u} \invDF D\stest{u} \invDF)\\
    &= \tr(-\invDF D\stest{u} \invDF D\test{u}).
\end{aligned}
\end{equation}
By making use of the rules above, we obtain the derivatives of the Lagrangian \cref{eq:Lagrangian} with respect to $u$ in the direction $\stest{u} \in V_0^p$:
\begin{equation}
\begin{aligned}
    \frac{\partial}{\partial u} L(u, \lambda) \stest{u} &= 
    \int_\Omega 
        (Du : Du)^{\frac{p-2}{2}} (Du : D\stest{u}) \; d x 
    + J'(\Omega)\, \stest{u} \\
    &+ ( \lambda_1, \dots, \lambda_d)^T \cdot \int_\Omega 
        \stest{u} \, \det(DF) 
        + (x + u) \tr(\invDF D\stest{u}) \det(DF) \; d x \\
    &+ \lambda_{d+1} \int_\Omega 
        \tr(\invDF D\stest{u}) \det(DF) \; d x.
\end{aligned}
\label{eq:Lu}
\end{equation}
Together with the usual derivative with respect to $\lambda$ into direction $\stest{\lambda} \in \R^{d+1}$, the optimality system reads
\begin{equation}
\begin{aligned}
    \frac{\partial}{\partial u} L (u, \lambda) \stest{u} &= 0  &&\quad \forall \; \stest{u} \in V_0^p\\
     \frac{\partial}{\partial \lambda} L (u, \lambda) \stest{\lambda} &= 0 &&\quad \forall \; \stest{\lambda} \in \R^{d+1}.
\end{aligned}
\label{eq:optimality_system}
\end{equation}
As the derivatives with respect to $\lambda$ can directly be taken form \cref{eq:Lagrangian}, we omit the details here.
In order to solve the nonlinear system \cref{eq:optimality_system} we require the linearization 
\begin{align}
    \frac{\partial^2}{\partial u^2} L(u^k, \lambda^k)(\stest{u}, \test{u}) &+  \frac{\partial}{\partial \lambda\,\partial u} L(u^k, \lambda^k)(\stest{u}, \test{\lambda}) &= - \frac{\partial}{\partial u} L(u^k, \lambda^k) \,\stest{u} && \forall \stest{u} \in V^p_0 \label{eq:LagrangianNewtonSystem1} \\
    \frac{\partial}{\partial u\,\partial \lambda} L (u^k, \lambda^k)(\test{u}, \stest{\lambda}) & &= - \frac{\partial}{\partial \lambda} L (u^k, \lambda^k) \,\stest{\lambda}  && \forall \stest{\lambda} \in \R^m \label{eq:LagrangianNewtonSystem2}
\end{align}
and the updates
\begin{equation}
    u^{k+1} = u^k + \test{u}, \quad \lambda^{k+1} = \lambda^k + \test{\lambda}
    \label{eq:LagrangianNewtonSystem3}
\end{equation}
where
\begin{equation}
\begin{aligned}
    &\frac{\partial^2}{\partial u^2}L (u,\lambda) (\test{u}, \stest{u}) = \\
    &\int_{\Omega} 
        (p-2)(Du : Du)^{\frac{p-4}{2}} (Du : D\test{u}) (Du : D\stest{u}) 
        + (Du : Du)^{\frac{p-2}{2}} (D\test{u} : D\stest{u}) \; d x \\
    &+ ( \lambda_1, \dots, \lambda_d)^T \cdot \int_\Omega \Big(
        \test{u} \, \tr(\invDF D\stest{u}) 
        + \stest{u} \, \tr(\invDF D\test{u}) \\
        &\qquad + (x + u) \left( 
            \tr(-\invDF D\stest{u} \invDF D\test{u})
             + \tr(\invDF D\test{u}) \tr(\invDF D\stest{u})
        \right) \Big) \det(DF) \; d x\\
    &+ \lambda_{d+1} \int_{\Omega} \Big( 
        \tr(-\invDF D\stest{u}\invDF D\test{u})
        + \tr(\invDF D\test{u})\tr(\invDF D\stest{u}) \Big) \det(DF) \; d x.
\end{aligned}
\label{eq:Luu}
\end{equation}
Reviewing the first integral in \cref{eq:Luu}, one observes that these terms do not exist for $p<4$, where $Du:Du=0$ holds on a set of non-zero measure.
However, this issue does not appear in the defect equation \cref{eq:Lu}, since there all exponents are non-negative.
We thus modify the first integral in \cref{eq:Luu} to
\begin{equation}
    \begin{aligned}
        \int_\Omega (p-2)( Du : Du + \epsilon \Theta(4 - p))^\frac{p-4}{2} (D\stest{u} : Du) (D\test{u} : Du)\\
        + (Du : Du + \epsilon )^\frac{p-2}{2} (D\stest{u} : D\test{u}) \; d x
    \end{aligned}
    \label{eq:modified-hessian}
\end{equation}
where $\Theta$ denotes the Heaviside function and $\epsilon>0$ a sufficiently small constant.
Notice that, within Newton's method in \cref{eq:LagrangianNewtonSystem1,eq:LagrangianNewtonSystem2}, this modification only affects the linearization and not the defect.
Thus, solutions of the original problem \cref{eq:DefRelaxedSteepestDescentDirection} are still obtained upon convergence.
Adding $\epsilon$ in \cref{eq:modified-hessian} serves two purposes. On the one hand guarantees invertibility and on the other it prevents divide-by-zero operations in the first term.

\section{Optimization Algorithm}
\label{sec:algorithm}

In this section we describe an algorithm for the solution of \cref{eq:DefRelaxedSteepestDescentDirection}. 
By the restriction of descent directions to maintain $g(F(\Omega))=0$, it is guaranteed that the geometric constraints are fulfilled up to a given tolerance at each iteration of the optimization process and not only on the limit.
The geometric constraints considered here, i.e. barycenter and volume of a free floating obstacle are particularly challenging to handle.
In an augmented Lagrangian or even pure penalty approach, the violation of $g=0$ in one iteration might lead to a strong overshoot of the shape deformation.
This causes oscillation of the shape because the geometry is unfeasible in each iteration.
For example, in particular at low Reynolds number flows, a major influence to the minimization of the energy dissipation is associated with the displacement of the flow by the obstacle.
Also minimizing the volume minimizes the energy dissipation.
At higher Reynolds number flows a descent direction is to move the obstacle downstream.
From a practical point of view this can only be solved by carefully adjusting initial values of the multipliers $\lambda$, the penalty factors, and the penalty increment values.
Thus, the practical attractiveness of the approach outlined here is that there are less heuristic and problem-dependent quantities to be adjusted. 
The user only has to provide the convergence criteria, the parameters of the step size control, and the values corresponding to the sequence of $p$, i.e. $p_\mathrm{max}$ and $p_\mathrm{inc}$.

From a mathematical point of view, the computational price one has to pay is the following: The set of admissible descent directions is not convex anymore, but the solution manifold of the nonlinear equation $g(F(\Omega))=0$.
For example, having computed an admissible step $u_p$ does not imply that $\tfrac{1}{2} u_p$ is also admissible.
This makes a step size control expensive, since the geodesics on the solution manifold are not straight lines in this case.
In \cref{alg:optimizationAlgorithm} the step size control is thus handled by scaling the shape sensitivity $J^\prime(\Omega)$ with a decreasing sequence $\sigma = (1, \tfrac{1}{2}, \tfrac{1}{4}, \dots)$.

As a note, for numerical reasons it might prove profitable to multiply equation \cref{eq:optimality_system} with $1/\sigma$.
\begin{algorithm}
\caption{Shape Optimization Steepest Descent Method}
\label{alg:optimizationAlgorithm}
\begin{algorithmic}[1]
    \Require{$\Omega$, $p_\mathrm{max}$}
    \State{$y \gets$ Solve primal problem}
    \State{$y_0 \gets y$}
    \State{Compute objective $\Phi_0 = J(\Omega)$}
    \Repeat
        \State{$y^\ast \gets$ Solve adj.\ problem}
        \State{$\sigma \gets 1$}
        \While{True}
            \State{$p \gets 2$}
            \State{$\bar u \gets 0$}
            \While{$p \le p_\mathrm{max}$}
                \State{($u_p,\lambda) \gets \textsc{NewtonSolver}(\bar u, \sigma, y, y^\ast)$}
                \State{$\bar u \gets u_p$}
                \State{Increase $p$}
            \EndWhile
            \State{Update geometry $\Omega$ with $u_{p_\mathrm{max}}$}
            \State{$y \gets$ Solve primal problem}
            \State{Compute objective $\Phi = J(\Omega)$}
            \If{$\Phi \geq \Phi_0$}
                \State{Update geometry $\Omega$ with  $-u_{p_\mathrm{max}}$}
                \State{$\sigma \gets \sigma / 2$}
                \State{$y \gets y_0$}
            \Else{}
                \State{$\Phi_0 \gets \Phi$}
                \State{$y_0 \gets y$}
                \State{\textbf{break}}
            \EndIf
        \EndWhile
    \Until{$\|u_{p_\mathrm{max}}\|_{W^{1,p}(\Omega)} < \epsilon_1$}
\end{algorithmic}
\end{algorithm}
\begin{algorithm}
\caption{Newton's Method for $p$-Laplacian Problem} \label{alg:NewtonSolver}
\begin{algorithmic}[1]
\Function{NewtonSolver}{$u_p, \sigma, y, y^\ast$}
	\State{$\lambda \gets 0$}
    \Repeat
        \State{($A$,$B$,$r_u$,$r_\lambda$) $\gets$ Assemble(u,$\lambda, y, y^\ast$) according to \eqref{eq:discreteDifferentialOperatiors}}
        \State{$(\test{u_p}, \test{\lambda}) \gets \textsc{SchurSolver}(A, B, r_u, r_\lambda, \test{u_p}, \test{\lambda})$}
        \State{$u_p \gets u_p + \test{u_p}$}
        \State{$\lambda \gets \lambda + \test{\lambda}$}
    \Until{$\| \test{u_p} \|_{W^{1,p}(\Omega)} + \| \test{\lambda} \|_2 < \epsilon_2$}
    \State{\textbf{return} ($u_p,\lambda$)}
\EndFunction
\end{algorithmic}
\end{algorithm}
\begin{algorithm}
\caption{Schur Complement Product} \label{alg:SchurComplementProduct}
\begin{algorithmic}[1]
\Function{SchurComplementProduct}{$A$, $B$, $w$}
    \For{$i = 1, \ldots, m$}
        \State{$b \gets b + B(:,i) w_i$}
    \EndFor
    \State{Solve $A z = b$}
    \For{$i = 1, \ldots, m$}
        \State{$b_i \gets - B(:, i)^T z$}
    \EndFor 
    \State{\textbf{return} $b$}
\EndFunction
\end{algorithmic}
\end{algorithm}
For the sake of readability we abbreviate the linearized optimality system \cref{eq:LagrangianNewtonSystem1,eq:LagrangianNewtonSystem2}  using the symbols
\begin{equation}
    \begin{aligned}
    A\test{u} &:= \frac{\partial^2}{\partial u^2} L(u^k, \lambda^k)(\stest{u}, \test{u}) &&\quad \forall \; \stest{u} \in V_{0,h}^{p}\\
    B\test{\lambda} &:= \frac{\partial}{\partial u\,d\lambda} L(u^k, \lambda^k)(\stest{u}, \test{\lambda}) &&\quad \forall \; \stest{u}\in V_{0,h}^{p}\\
    B^T\test{u} &:=\frac{\partial}{\partial \lambda\,\partial u} L(u^k, \lambda^k)(\test{u}, \stest{\lambda}) &&\quad \forall \; \stest{\lambda}\in \R^d\\
    r_u &:= - \frac{\partial}{\partial u} L(u^k,\lambda^k) \,\stest{u} &&\quad \forall \; \stest{u}\in V_{0,h}^{p}\\
    r_\lambda &:= -\frac{\partial}{\partial \lambda} L(u^k,\lambda^k) \,\stest{\lambda} &&\quad \forall \; \stest{\lambda} \in \R^d.
    \end{aligned}
	\label{eq:discreteDifferentialOperatiors}
\end{equation}
With $V_{0,h}^{p}$ a discrete approximation of $V_{0}^{p}$ used for a finite element discretization of \cref{eq:LagrangianNewtonSystem1,eq:LagrangianNewtonSystem2}, which then leads to the saddle point problem.
\begin{equation}
    \begin{pmatrix}
    A & B \\
    B^T & 0
    \end{pmatrix}
    \begin{pmatrix}
    \test{u} \\
    \test{\lambda} 
    \end{pmatrix} = 
    \begin{pmatrix}
    r_u \\
    r_\lambda
    \end{pmatrix}
    \label{eq:DiscreteSaddlePointProblem}
\end{equation}
where $A\in \R^{n\times n}$ and $B\in \R^{n\times m}$. 
In order to solve for the increments $\test{u}$ and $\test{\lambda}$, we formally apply one block wise Gauss elimination and obtain
\begin{equation}
    \begin{pmatrix}
    A & B \\
    0 & -B^T A^{-1}B
    \end{pmatrix}
    \begin{pmatrix}
    \test{u} \\
    \test{\lambda} 
    \end{pmatrix} = 
    \begin{pmatrix}
    r_u \\
    r_\lambda - B^T A^{-1} r_u
    \end{pmatrix}
    \label{eq:SchurComplementReducedSystem}
\end{equation}
where $S:=-B^T A^{-1}B$ is the so-called Schur complement operator.
In order not to explicitly compute $A^{-1}$ a equation system with $A$ is solved instead.
In general, the optimality system \cref{eq:optimality_system} of problem \cref{eq:generalShapeOptProp1,eq:generalShapeOptProp2,eq:generalShapeOptProp3} is highly nonlinear.
Especially with increasing values of $p$ the solution process becomes more challenging unless a good initial guess $u_p^0$ is provided.
To overcome this issue, and to reduce computational effort, we consider a finite sequence $p_k := p_\mathrm{init} + k p_\mathrm{inc}$ where $p_\mathrm{init} := 2$. First, the solution for $p_\mathrm{init}$ with initial $u_{p_\mathrm{init}} = 0$ and $\lambda = 0$ is computed. Thereafter the solution of the constraint $p_k$-Laplacian problem \cref{eq:DefRelaxedSteepestDescentDirection} is used as an initial guess for the $p_{k+1}$-Laplacian problem, cf. \cite{mueller2021}.
Here, with the choice of $p_\mathrm{max}$ we adjust the approximation quality of Lipschitz deformations.

The overall optimization procedure is outlined in \cref{alg:optimizationAlgorithm}. 
The steepest descent method is reflected in the loop spanning from lines 4 to 29, where the necessary optimality condition is checked.
Here, $y$ again denotes the state variable of the PDE constraint $e$, which we refer to as the primal problem.
Nested within this loop, a step-size control operates in the lines 8 to 28.
It checks whether the proposed next shape $F(\Omega)$ leads to an improvement of the objective function, in terms of the displacement field $u_p$.
If not, then the parameter $\sigma$ is reduced.
Note that, in contrast to the classical backtracking line search in linear spaces, we have to recompute the descent direction $u_p$ hereafter.
This is due to the fact that by shortening the step-length, we can not follow straight lines towards $0 \in V_0^p$, but have to stay within the solution manifold of the non-linear geometric constraints $g(F(\Omega))=0$.
In line 5 the adjoint PDE is solved, which yields the adjoint state $y^\ast$.
After this, the shape sensitivity $J^\prime(\Omega)$ can be evaluated in line 6, which depends on $y$ and $y^\ast$.

In line 11 of \cref{alg:optimizationAlgorithm}, the nonlinear solver for the steepest descent problem \cref{eq:DefRelaxedSteepestDescentDirection} is called, which can be seen in \cref{alg:NewtonSolver}.
The key part of this solver is the solution to the saddle point problem \cref{eq:DiscreteSaddlePointProblem} in the Schur complement form \cref{eq:SchurComplementReducedSystem}.
This could be realized with a variety of iterative solvers, which are not further specified here. Popular approaches for these kind of problems are the Uzawa iteration and the Arrow-Hurwicz algorithm.
For this purpose, \cref{alg:SchurComplementProduct} outlines the computational steps for a matrix-vector product with $S$.

\section{Numerical Methodology}
Results of the present study are obtained from the open-source toolbox UG4~\cite{Vogel2014ug4}. 
This simulation framework has MPI-based parallelization, and features a geometrical multigrid preconditioner~\cite{Hackbusch85}.
The grid partitioning and load balancing scheme is based on ParMetis~\cite{PARMETIS}.

Stable $P_2-P_1$ finite elements were used to discretize the governing nonlinear Navier-Stokes equations \cref{eq:NavStokes} and their linearization, therefore no additional stabilization is required. Moreover, the viscosity is $\nu=0.02$ in all cases.
The same setting was used to discretize the linear adjoint problem, cf. ~\cite{pinzon2021,onyshkevych2020} and the references therein.
As regards the p-Laplace relaxation problem, for which the optimality system is described in \cref{eq:SchurComplementReducedSystem}, $P_1$ Lagrange shape functions were employed.
Computational grids consist of triangular (2d) and tetrahedral elements (3d). 
They were generated using GMSH~\cite{GMSH}. 
 
The simulations followed the workflow proposed in \cref{alg:optimizationAlgorithm}.
At the beginning of each optimization step, the steady, incompressible flow, described by Navier-Stokes equations \cref{eq:NavStokes}, was computed followed by the solution of the corresponding adjoint system.
The p-Laplacian descent algorithm  \cref{alg:optimizationAlgorithm} initially employed  $p_\mathrm{init}=2.0$ and incremented $p$ by $p_\mathrm{inc}=0.19$. The given maximum values of $p$ read $p_\mathrm{max}=4.8$ [$p_\mathrm{max}=4.1$] for the computed 2d [3d] test cases. 
Termination criteria of \cref{alg:optimizationAlgorithm} and \cref{alg:NewtonSolver} were always set as $\epsilon_1=\num{1E-5}$ and $\epsilon_2=\num{1E-8}$, respectively. 
The modification term introduced in \cref{eq:modified-hessian} reads  $\epsilon=\num{1E-8}$ for all cases. 

As a practical note, care must be taken to correctly interpolate the values of $\vel$, $\press$, and their respective adjoints.
These are involved in the assembly of $J'$, which is present in \cref{eq:SchurComplementReducedSystem}.
The geometrical constraints are part of this system of equations.
It has been described that they lead to an $m\times m$ system of equations, so their discretization is not within a finite element space but in $\R^m$.
For the investigated case cases, $m\in \{3,4\}$ in 2d and 3d, respectively. 
Thus, we use a direct solver to find the solution of the Schur complement system.

The corresponding codes used for these results can be found in the online repository \cite{gmgshapeopt2021}.

\section{Results}
\label{sec:numerics}
This section presents results for 2d and 3d fluid dynamics applications. They either refer to an initial square (2d) or cube (3d) centrally placed in a rectangular flow domain at low Reynolds number, i.e. $Re=1\cdot H/\nu$ = 20, where $H$ refers to the length of the initial edges. The employed box-domain is outlined in Fig. \ref{fig:domain}. It spans 20 units in length and $\gamma=6$ units in height (2d, 3d) and depth (3d), respectively, and the flow enters the domain through the left vertical boundary. 
The inflow profile on $\Gin$ features a peak unit-value in the center of the inlet plane and is described by
\begin{equation*}
	\vel_\infty=\left(\max\left\{0, \prod\limits_{i=2}^d \cos(\frac{\pi|x_i|}{\delta})\right\}, 0, \dots, 0\right) \in \R^d
\end{equation*}
where $\delta$ corresponds to the inlet height.

The central aspect of the paper is the creation and removal of geometrical singularities.
Emphasis is placed on illustrating and explaining how the corners of the obstacle are removed during the optimization process, as well as how tips are generated to reach an optimal shape. 
%
A crucial aspect is the evolution of the mesh quality during an optimization. We utilize the 2d studies, to compare the mesh quality of the optimal and the initial design by means of the ratio $\rho$ between the radii of circumcircle and incircle, and report the extreme values of the interior angles of the triangulation.
Moreover, we describe the behavior of the proposed algorithm in two different 3d configurations, where the surface of the obstacle is highly resolved.


Mind that the geometrical constraints are preserved during each optimization step for every value of $p$, since they are incorporated to the system of equations. 
Their fulfillment is included in the convergence condition set for Newton's method in \cref{alg:NewtonSolver}, therefore there is no need to provide results for their fulfillment per step.
Solving the nonlinear system of \cref{eq:optimality_system} implies solving the geometrical constraints \cref{eq:BarycenterConstraint,eq:VolumeConstraint} to the error reduction tolerance set for Newton's method.
%
%
The major portion of the computational effort in \cref{alg:optimizationAlgorithm} is spent on solving the  $p$-Laplace relaxed problem via the scheme described in \cref{sec:algorithm}.  
Particularly, lines 11-15 of \cref{alg:optimizationAlgorithm} are computationally expensive, as will be explained here and in \cref{sec:scalability}.

  
\begin{figure}[!htbp]
    \center\includegraphics[width=0.8\textwidth]{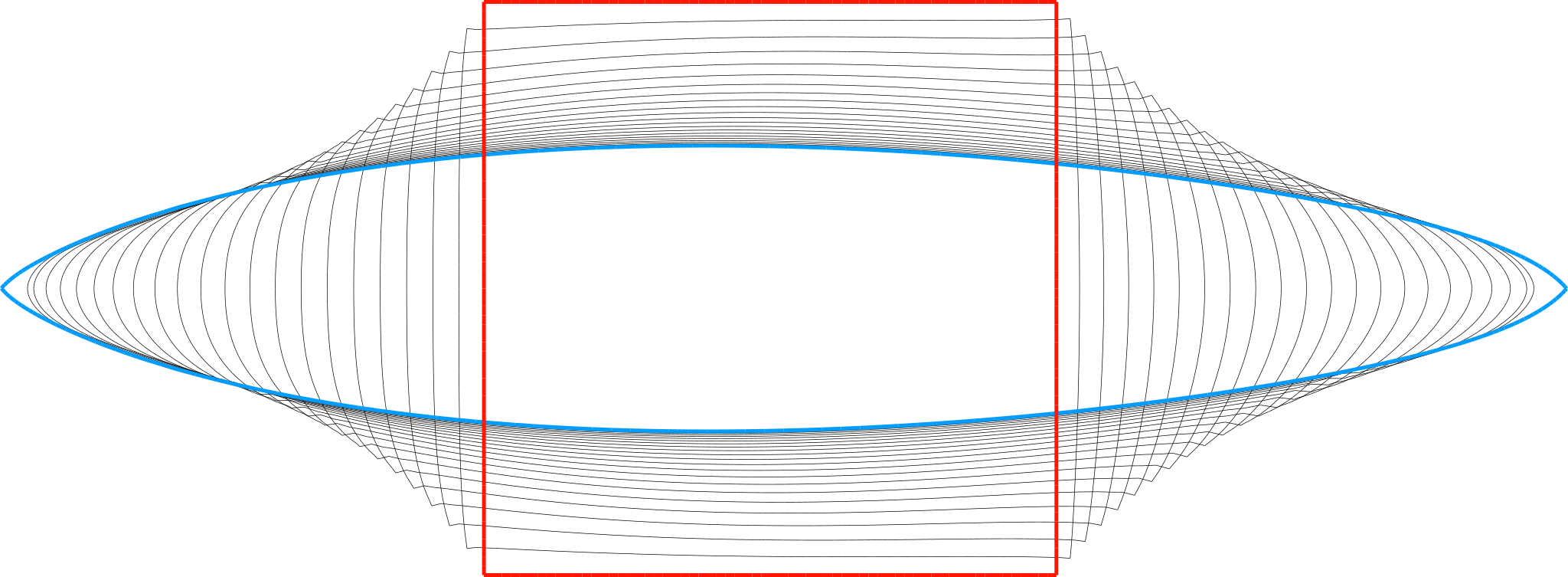}
    \caption{Superposition of the deformation sequence for a 2d configuration. The obstacle's initial shape {\color{red} red} is presented superimposed to the sequence of generated shapes {\color{gray}gray} until an optimal shape {\color{blue}blue} is obtained upon convergence.}
    \label{fig:2d-morphing}
\end{figure}

\subsection{Two-dimensional studies}
\begin{figure}[!htbp]
    \center\includegraphics[width=0.8\textwidth]{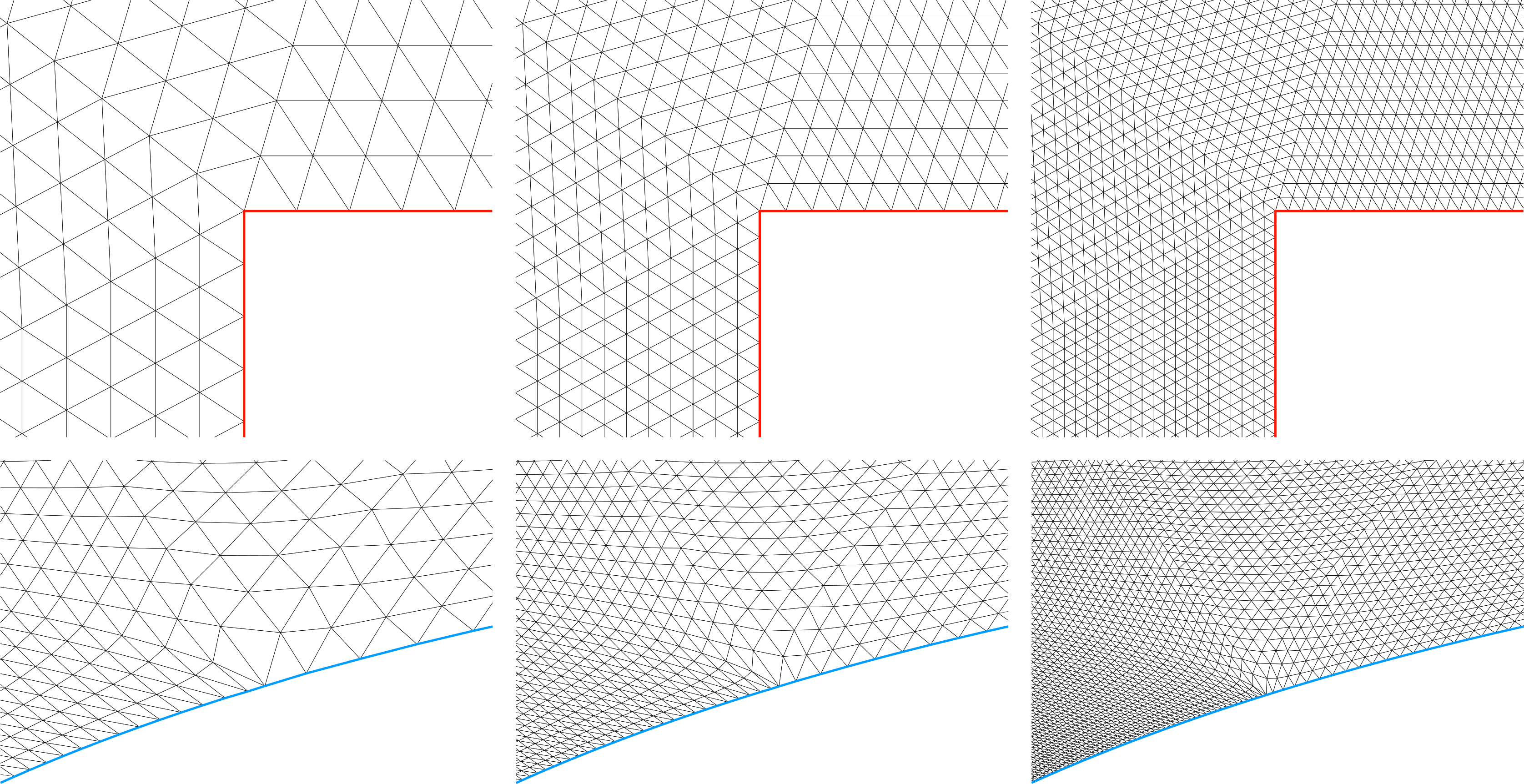}
    \caption{Removal of the geometrical singularity in the obstacle's initial configuration across several levels of refinement. For 4, 5, and 6 refinements the upper left corner of the box is smoothed via updating the geometry $\Omega$ iteratively, as stated in line 16 in \cref{alg:optimizationAlgorithm}.}
    \label{fig:2d-corners}
\end{figure}

Simulations in the 2d domain were performed for several levels of refinement, to better describe the removal of the geometrical singularities, as well as the mesh quality. Figure \ref{fig:2d-morphing} compares the 
 initial design (red) with the converged design (blue), together with a contour plot of a deformation sequence (gray).
A robust removal of the box corners is clearly visible, as well as the creation of the tips in the rear and the aft sections.
As described in \cref{sec:model}, the geometrical constraints are preserved in all optimization steps.
This feature can be observed by the continuous transition between shape iterates until an optimum is obtained.
In contrast, in~\cite{pinzon2021} bouncing of the shapes during the early stages of the optimization is reported, which is related to an approximate solution of the geometrical constraints.
Figure \ref{fig:2d-corners} magnifies, the geometry and the mesh in the upper-left corner of the obstacle. 
The initial and final shapes are presented on the top and bottom, respectively, for different grid refinement levels from left to right. 
The figure displays that the smoothing occurs similarly on all grids, and the elements around the initial singularity are not dramatically degenerated during the optimization. 
Towards the last step, no indication of the initial geometric singularity is visible on the obstacle's surface.  

The mesh quality is investigated for the final step using  4, 5, and 6 levels of refinement. 
As outlined in~\cref{sec:algorithm}, a series of shape iterates are obtained until an optimum, with respect to \cref{eq:EnergyDispObjectiveFunctional}, is found. 
Mind that the geometric multigrid preconditioner, which is used to allow for numerical scalability, requires the generation of a grid hierarchy, of which we provide the base level, i.e. the coarsest mesh. 
This implies that the simulations are based upon a predetermined mesh quality, and while the optimization we propose in~\cref{alg:optimizationAlgorithm} aims at preserving grid quality, it doesn't contemplate improving it with respect to the initial geometry. 
Table \ref{table:qualities} provides quality measurements for the final step, when the optimal shape is found, using several grid refinement levels. Assessed data refers to the worst triangular elements extracted from the 2d grid, i.e. the observed minimum and maximum interior angles, and the largest radius ratio.
 We also compare the radius ratio between the last and first configurations. 
The value of $\rho_0=1.468$ indicates that the initial mesh does not have an ideal quality. 
Results also demonstrate that, if $p$ is high enough, the approximation of Lipschitz transformations, as seen in \cref{eq:V-0-p,eq:V-0-infty}, prevents a significant loss of mesh quality over mesh refinements. For the presented 2d cases, a value of $p=4.8$ yielded a sufficient approximation to $p=\infty$ in terms of the mesh quality, while allowing for the creation and removal of geometrical singularities. 
The mesh refinement study might reveal that a higher maximum p-values are necessary for the finer grids, since the quality slightly deteriorates.
Nevertheless, numerical stability must be taken into account when increasing this value, given that it is used in \cref{eq:Luu} as an exponent.
The latter fact turned out to be a limiting factor in our numerical simulations.
However, the measurement of the worst minimum and maximum angles express that the triangles, which have undergone the largest deformation, are still not close to being critical.

\begin{table}
	\centering
	\begin{tabular}[t]{r r r r r r}
			Refinements & Elements      & Minimum angle & Maximum angle & Radius ratio $\rho_\infty$ & $\frac{\rho_\infty}{\rho_0}$\\ \hline
			4           & \num{70656}   & \num{13.41}   & \num{132.32}  & \num{3.20}   & \num{2.18}\\
			5           & \num{282624}  & \num{11.93}   & \num{139.03}  & \num{4.24}   & 
			\num{2.89}\\
			6           & \num{1130496} & \num{9.94}    & \num{145.04}  & \num{5.76}   & \num{3.92}\\
	\end{tabular}
	\caption{Assessment of mesh quality evolution for several refinement levels observed in 2d; 
	Displayed data for minimum and maximum interior angles supplemented by the largest radius ratio of the triangulation extracted for the last optimization step, where an optimal shape is reached. The last column compares the largest radius ratio of the 
	 optimal shape (subscript $\infty$) 
	 and the initial configuration (subscript $0$, $\rho_0 = \num{1.468}$).
	}
	\label{table:qualities}
\end{table}

\subsection{Three-dimensional studies}
\label{subsec:threeDResults}
Results for the 3d simulations refer to 4 levels of grid refinement. The computational grid has a total of \num{4980736} tetrahedrons, and \num{49152} triangles discretize the surface of the obstacle, $\Gobs$, on the highest refinement level.
Our optimization scheme generates a series of deformation fields $u_p$ that, applied to the domain $\Omega$, results in an optimal shape with respect to the energy dissipation \cref{eq:EnergyDispObjectiveFunctional}. 

\begin{figure}[!htbp]
	\begin{center}
		\begin{tabular}{ccc}
			(a) & \includegraphics[width=0.45\textwidth]{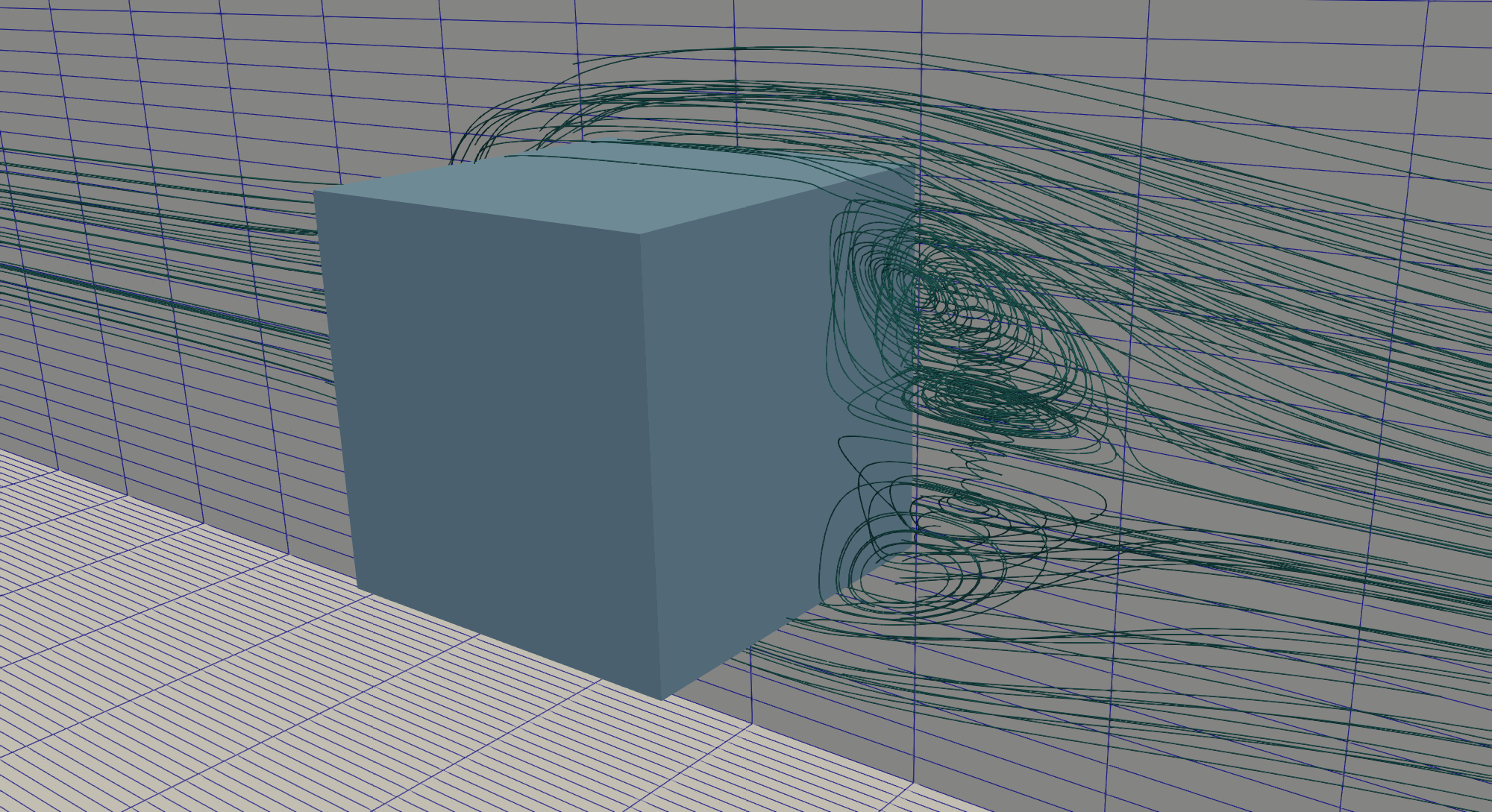}& \includegraphics[width=0.45\textwidth]{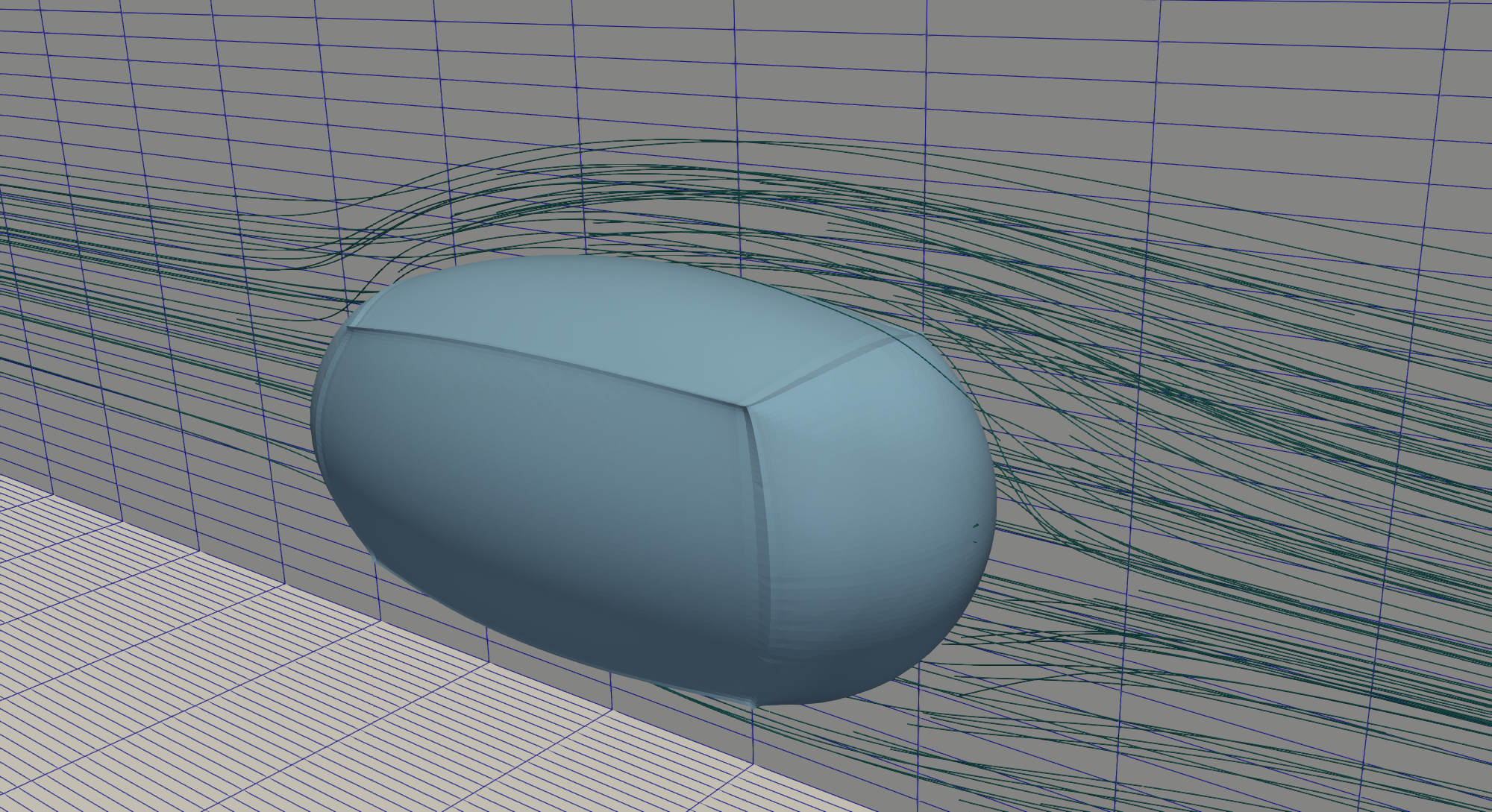}\\
			(b) &\includegraphics[width=0.45\textwidth]{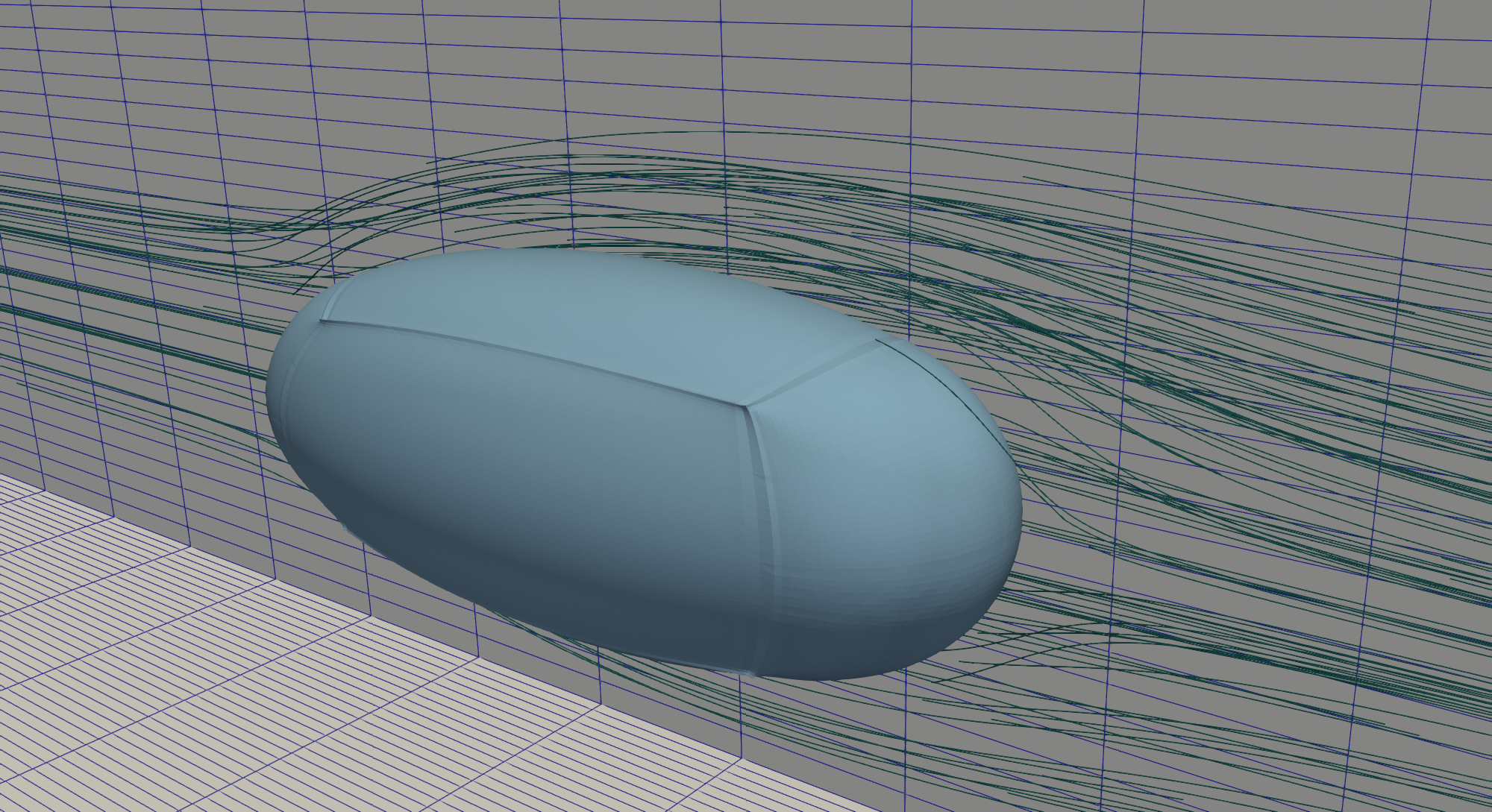}& \includegraphics[width=0.45\textwidth]{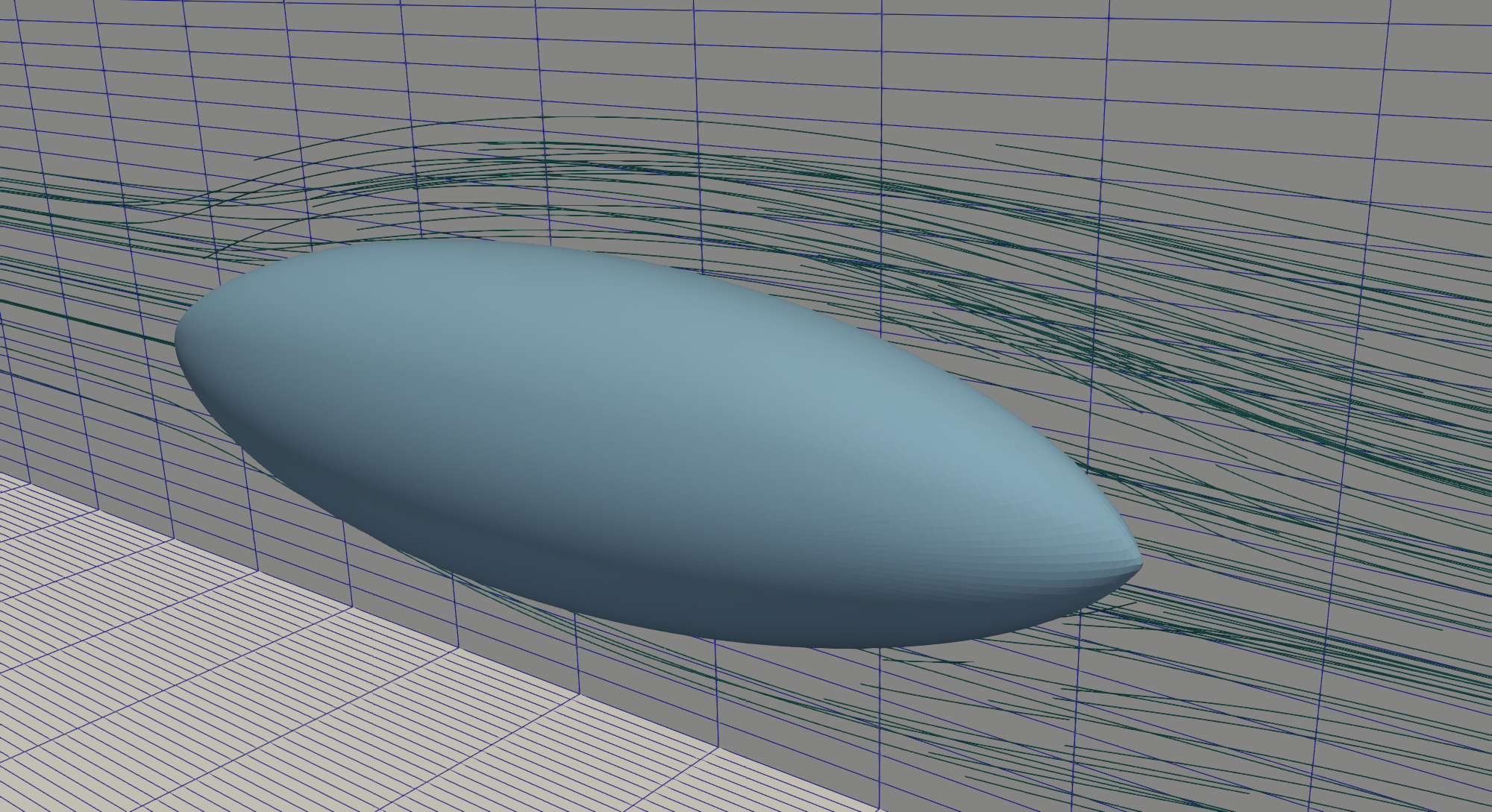}\\
		\end{tabular}
	\end{center}
	\caption{Streamlines for the rear of the obstacle located in the wind tunnel.}
	\label{fig:3d-streamlines}
\end{figure}

\Cref{fig:3d-streamlines} presents iterated shapes from the initial to the final optimization step.
It shows the downstream part of the geometry.
For the reference shape the streamlines visualize a region where the flow direction points backwards w.r.t.~the main flow direction.
Since this effect contributes to the energy dissipation it vanishes during the optimization at an early stage.
This phenomenon can be quantified by observing the shear stress acting on the surface of the obstacle $\tau \cdot e_1 = \left(\nu (Dv + Dv^T) \cdot n\right) \cdot e_1$.
Here $e_1$ is the first unit vector describing the main flow direction.
For the 3d (cf.~\cref{subsec:threeDResults}) case $\tau \cdot e_1 \in [-11.06, 1.27]$ for the initial shape and $\tau \cdot e_1 \in [-20.53, -1.12]$ for the final one, respectively.

\begin{figure}[!htbp]
	\begin{center}
		\begin{tabular}{lcccc}
			Step & Obstacle& Tip & Corners\\
			0 &\includegraphics[width=0.35\textwidth]{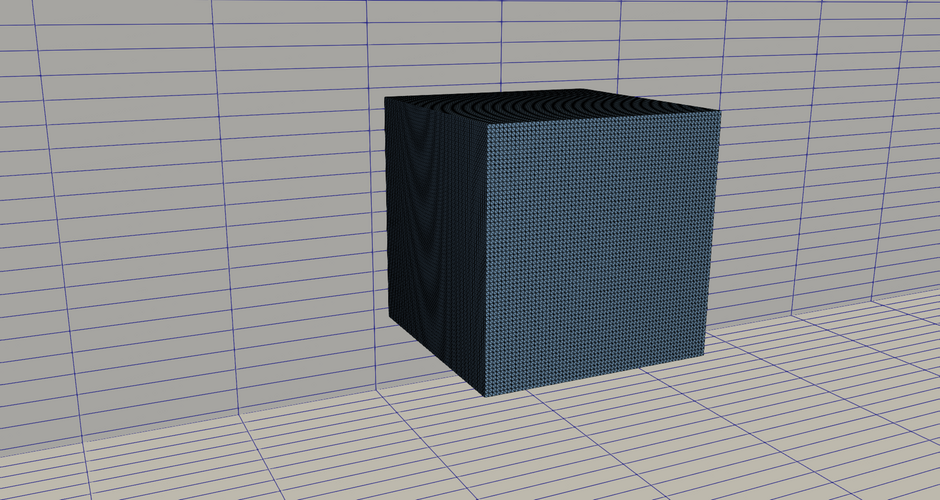}& \includegraphics[width=0.25\textwidth]{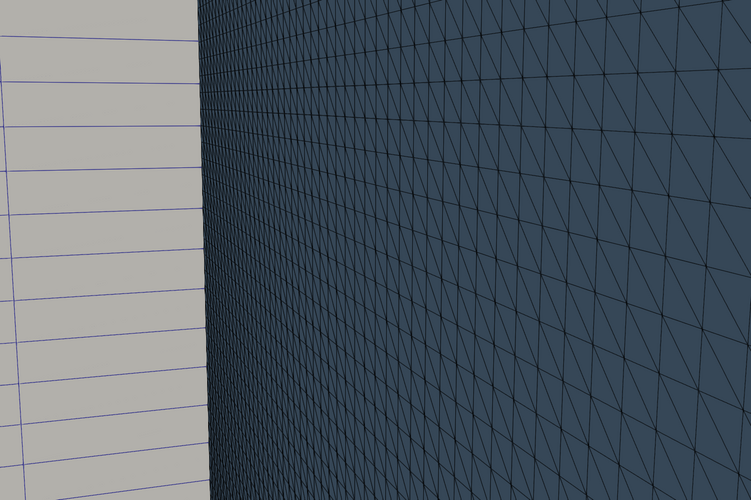}& \includegraphics[width=0.25\textwidth]{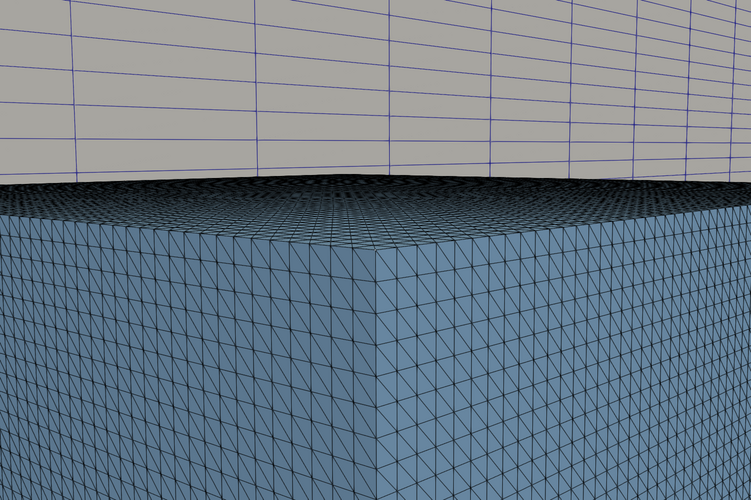}\\
			5 &\includegraphics[width=0.35\textwidth]{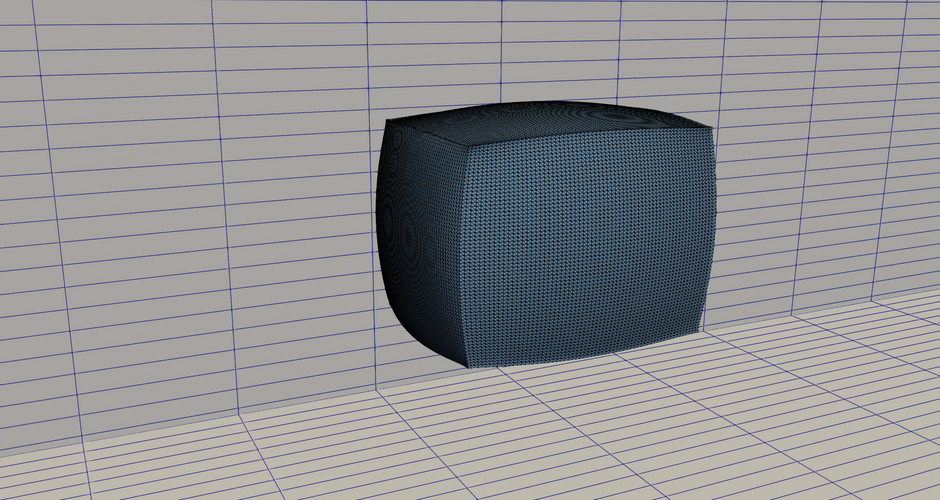}& \includegraphics[width=0.25\textwidth]{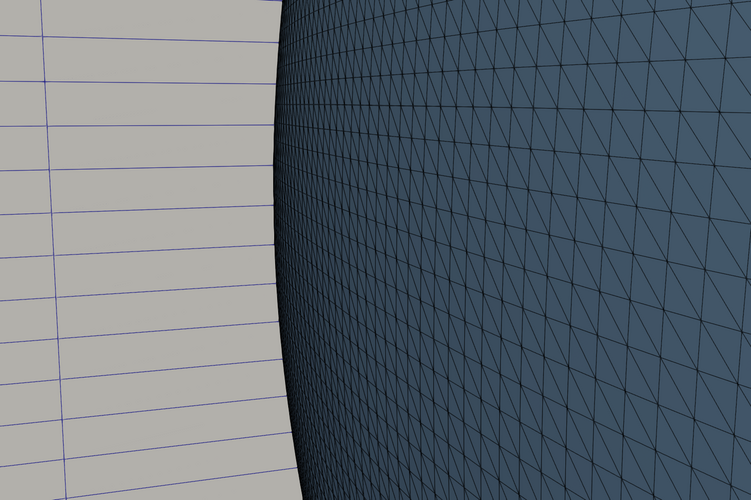}& \includegraphics[width=0.25\textwidth]{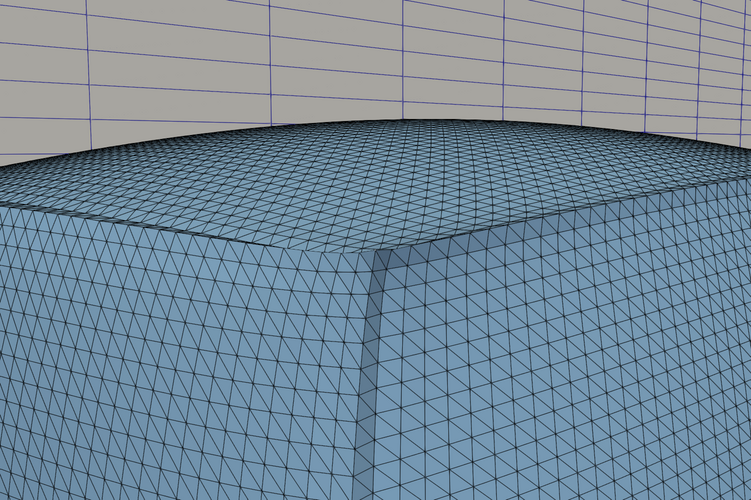}\\
			15 &\includegraphics[width=0.35\textwidth]{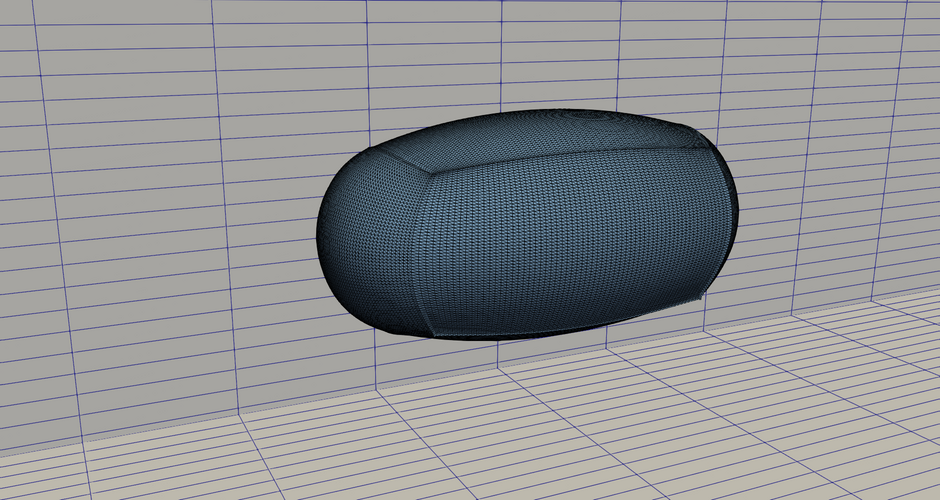}& \includegraphics[width=0.25\textwidth]{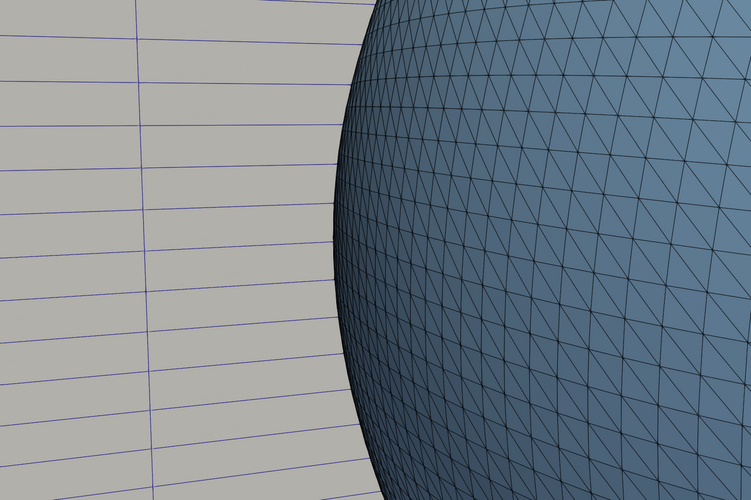}& \includegraphics[width=0.25\textwidth]{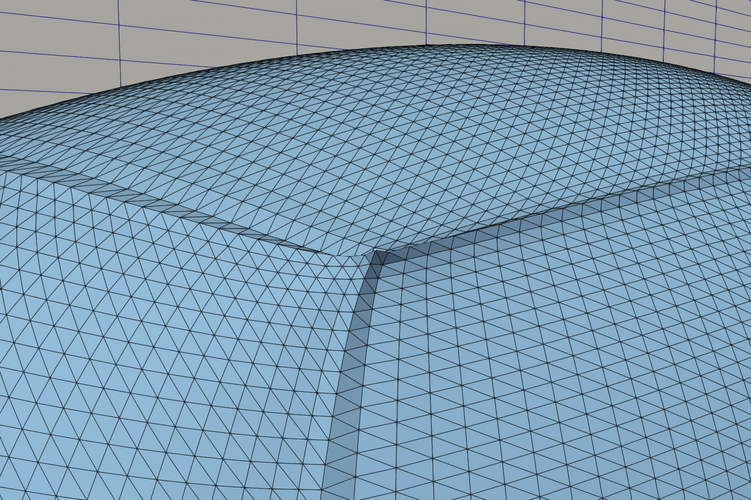}\\
			35 &\includegraphics[width=0.35\textwidth]{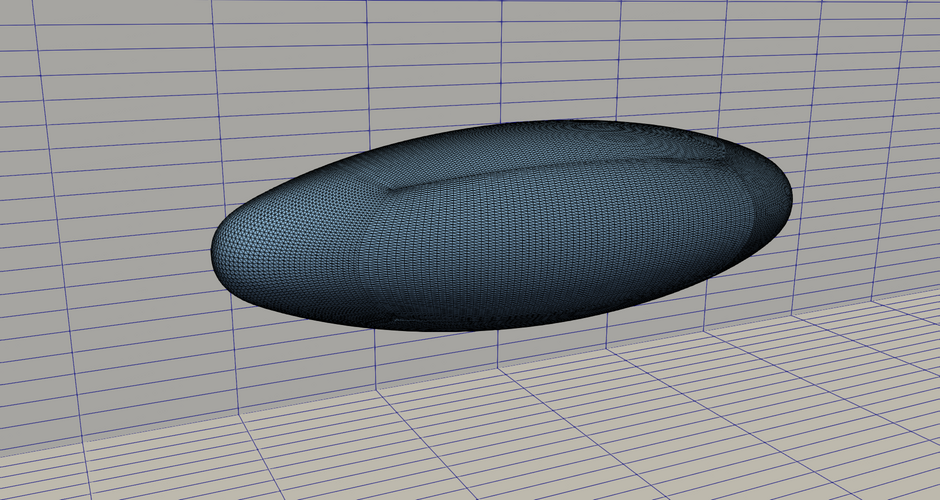}& \includegraphics[width=0.25\textwidth]{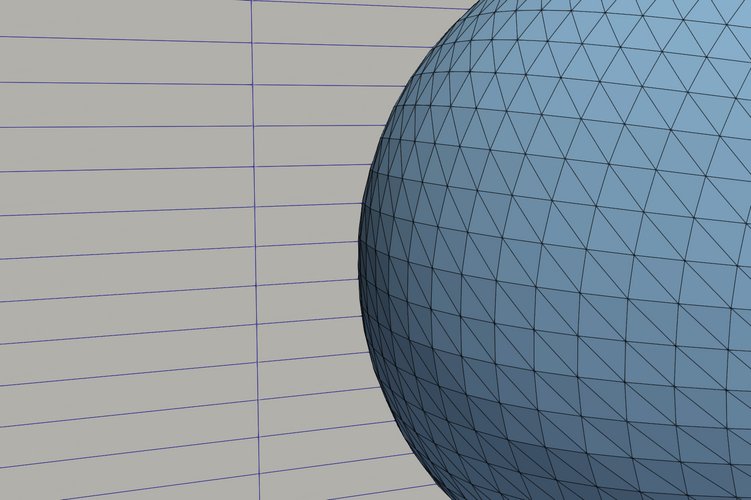}& \includegraphics[width=0.25\textwidth]{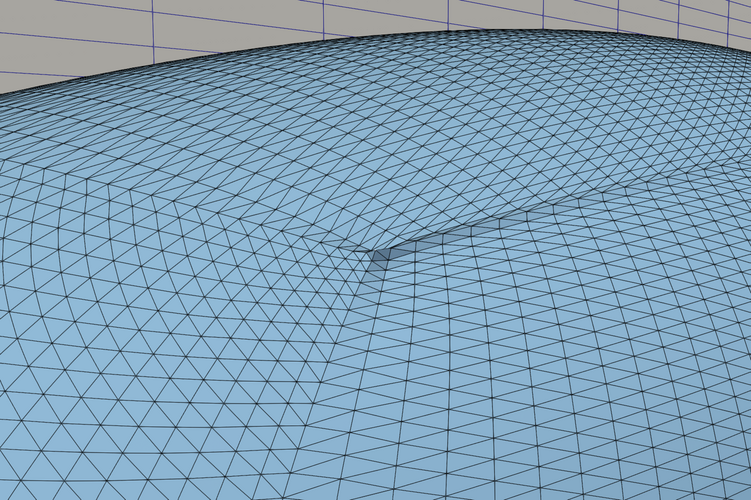}\\
			50 &\includegraphics[width=0.35\textwidth]{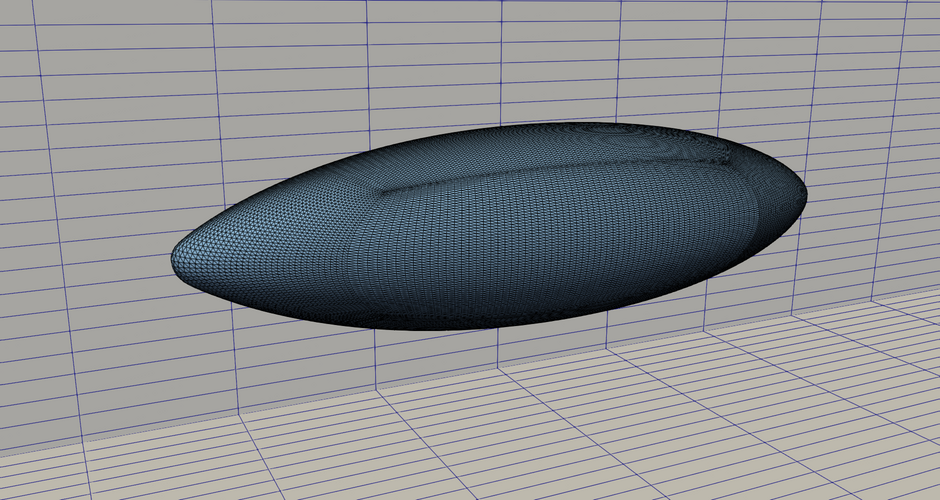}& \includegraphics[width=0.25\textwidth]{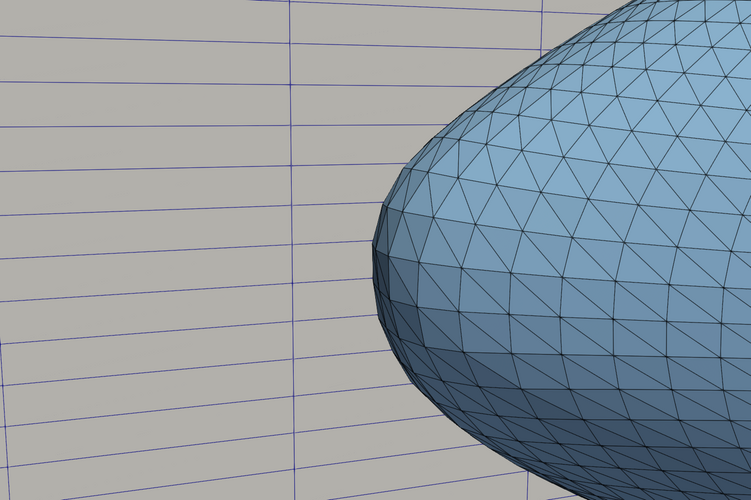}& \includegraphics[width=0.25\textwidth]{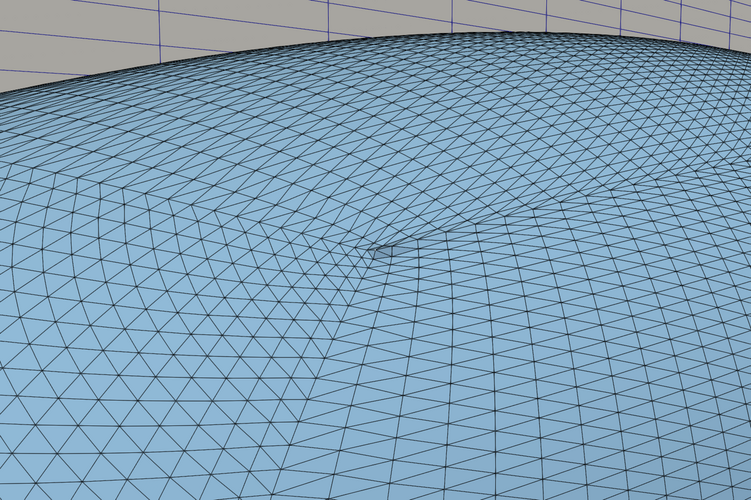}\\
			100 &\includegraphics[width=0.35\textwidth]{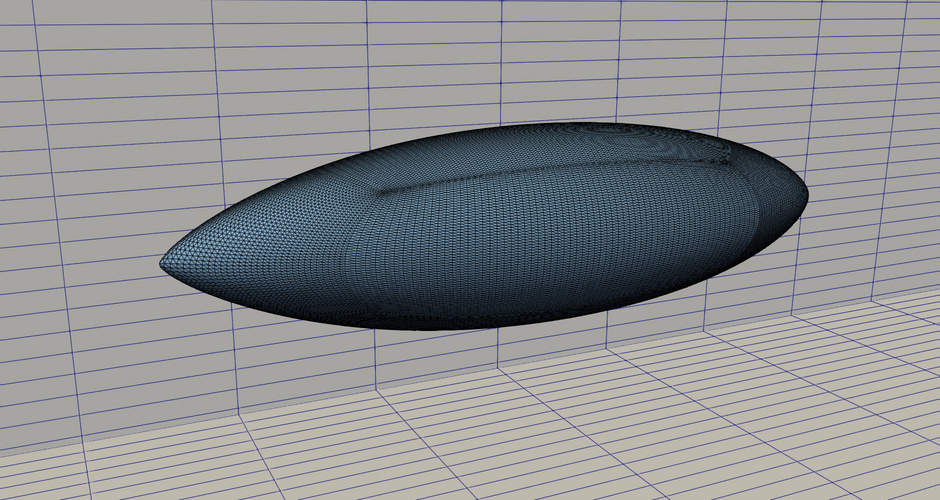}& \includegraphics[width=0.25\textwidth]{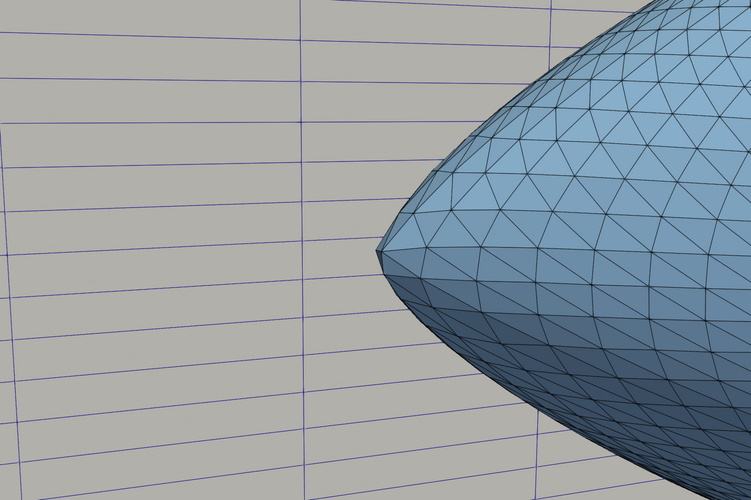}& \includegraphics[width=0.25\textwidth]{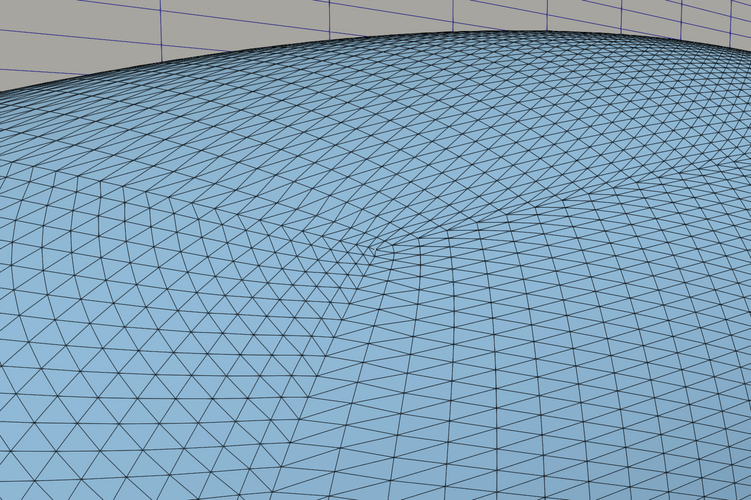}\\
			
		\end{tabular}
	\end{center}
	\caption{Deformation sequence for optimization steps \{0, 5, 15, 35, 50, 100\}. The complete obstacle, together with a detailed view of the geometrical singularity removal and generation process, are presented.}
	\label{fig:3d-sequence}
\end{figure}

Similar to the 2d case, the edges and corners displayed by the initial geometry are gradually removed as part of the optimization process, cf. \cref{fig:3d-sequence}. 
Additionally, the tips created at the central upstream and downstream ends shapes a streamlined body that does not feature any separation, as shown in row (b) of \cref{fig:3d-results-table}.  

\Cref{fig:3d-sequence} shows a deformation sequence of the 3d case, starting from the initial configuration and ending with an optimum obstacle surface.
The figure focuses on the overall shape (left), an exemplary corner of the initial geometry (center), as well as the location of the upstream end of the final geometry (right). 
During the initial steps, the obstacle aligns to the flow, i.e. is stretched in the direction of the flow and compressed in the other two directions. 
Edges begin to emerge from the round upstream and downstream facing surfaces, thus creating the geometry observed in step 15. 
Subsequently, a round cross section starts to take form in the center and as seen in step 35, where the final tip locations also become more apparent. 
Recall that the mesh deformation corresponds to line 15 in \cref{alg:optimizationAlgorithm}.
We again emphasize that 
all shape iterates meet the volume and barycenter constraints, which are deemed crucial for the success of this optimization scheme. Footprints of initial corners and edges are still visible in the mesh at later stages of the optimization, 
e.g. step 50.
However, they are completely smoothed out towards the end of the simulation and only the macro elements, resulting from the grid hierarchy, are visible.

The front tip is shown for step 100, where also the previously existing singularities have disappeared. 

\begin{figure}[!htbp]
	\centering
	\includegraphics[width=0.8\textwidth]{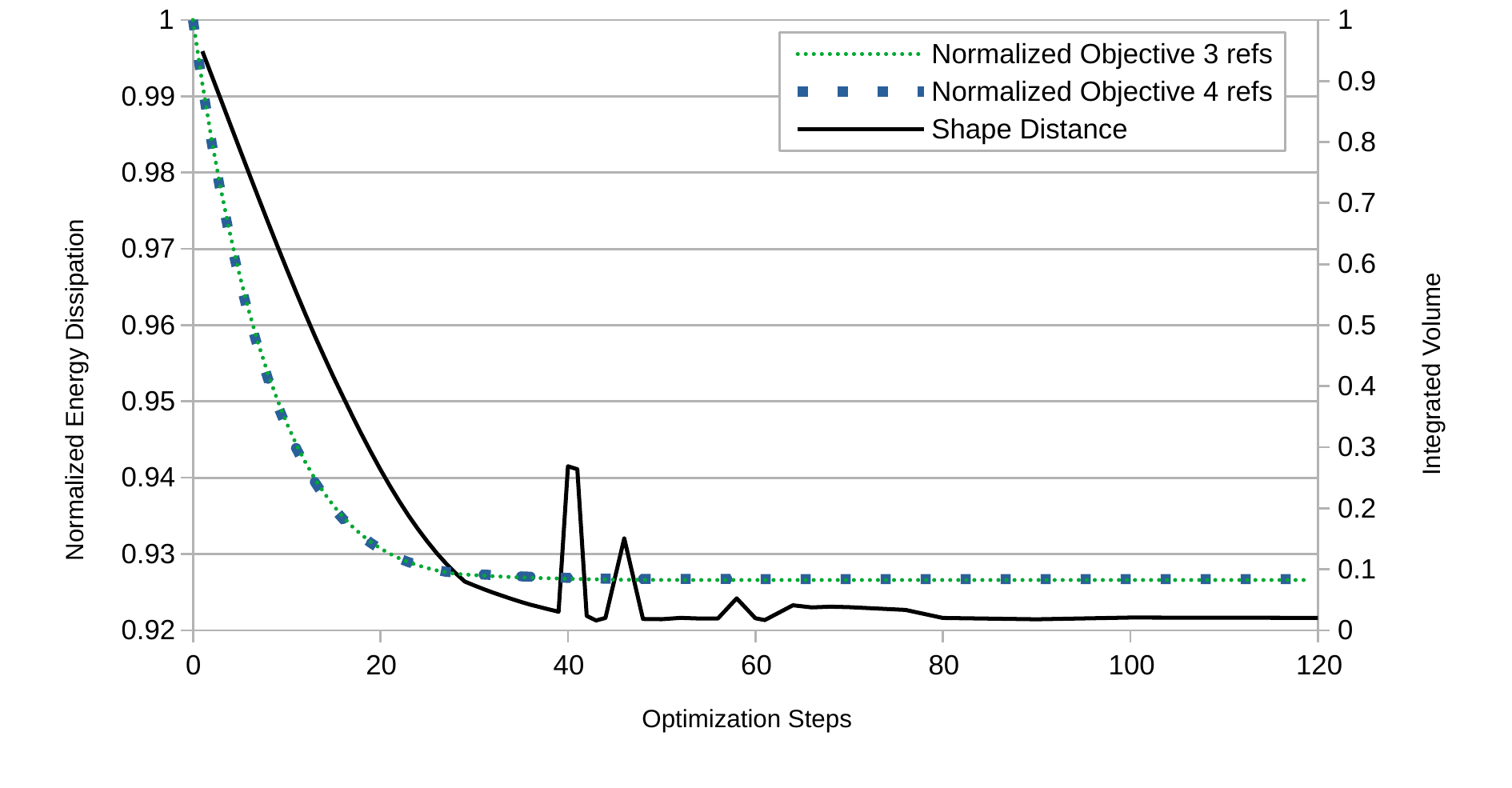}
	\caption{3d results for 3 and 4 levels of refinement are compared. The energy dissipation, see \cref{eq:EnergyDispObjectiveFunctional}, is plotted against the difference between shapes of each refinement level per step.}
	\label{fig:drag-distance}
\end{figure}

As mentioned in \cref{sec:model}, this approach optimizes the obstacle's shape for the functional given in \cref{eq:EnergyDispObjectiveFunctional}. 
Therefore, results are provided in \cref{fig:drag-distance} for a 3d setting, which show how the generated shape, after convergence of \cref{alg:optimizationAlgorithm}, consists of an optimum with respect to the cost function.

Figure \ref{fig:drag-distance} depicts the objective function plot evolution over 120 optimization steps using 3 and 4 levels of grid refinement, respectively.
The fact that the objective function \cref{eq:EnergyDispObjectiveFunctional} decreases monotonically is linked to lines 15-26 of \cref{alg:optimizationAlgorithm}, where a line search strategy is implemented.
Once the deformation field is obtained for $p_{max}$, the geometry is updated and we get a new obstacle shape. 
The state equation is solved and the cost function calculated to guarantee that the new shape iterate represents a descent direction.
As seen in lines 18-21, whenever the condition is false, the deformation is withdrawn and the step size control value is reduced to repeat the unsuccessful step with a scaled shape sensitivity $J'$.
As a further indicator for the convergence, we approximate the distance between the iterated shapes $\Omega^k_\mathrm{3ref}$ of the 3 refinements run to the optimal solution of the 4 refinements run $\Omega^\infty_\mathrm{4ref}$.
\Cref{fig:drag-distance} shows the integrated volume that refers to the symmetric difference
as 
\begin{equation}
	d(\Omega^k_\mathrm{3ref}, \Omega^\infty_\mathrm{4ref}) := \left|\Omega^k_\mathrm{3ref} \setminus  \Omega^\infty_\mathrm{4ref}\right| + \left|\Omega^\infty_\mathrm{4ref} \setminus  \Omega^k_\mathrm{3ref}\right|
\end{equation}
The integration is carried out using the boolean filters of the VTK library~\cite{VTK} with which a triangulation of the surface of the volumes of interest can be obtained.
A concatenation of the VTK boolean filters gives us a surface triangulation together with its normal vector.
Then, utilizing divergence's theorem the volume can be found.

\subsection{Scalability Study}
\label{sec:scalability}
Weak scalability of the solution strategy for the $p$-Laplacian relaxed problem, from $p_\mathrm{init}=2.0$ up to $p_\mathrm{max}$ is presented here.
This solution strategy, described in \cref{sec:algorithm}, is referred in these results as the $p$-solver.
This corresponds to lines 11-15 of \cref{alg:optimizationAlgorithm}. 
It was studied for up to \num{262144} cores in a 3d setting.
The study was carried out with the supercomputer Hawk at HLRS.
It features 5632 compute nodes, each with a dual-socket architecture and a total of 128 cores. Each core with a maximum frequency of 2.25GHz, and 256GB of RAM.
The runs were carried out taking into account the hypercube topology of the system to maximize core usage and minimize parallel communication overhead.

\begin{figure}[!htbp]
	\centering
	\begin{subfigure}[t]{\textwidth}
		\centering
 	    {\includegraphics[width=0.4999\textwidth]{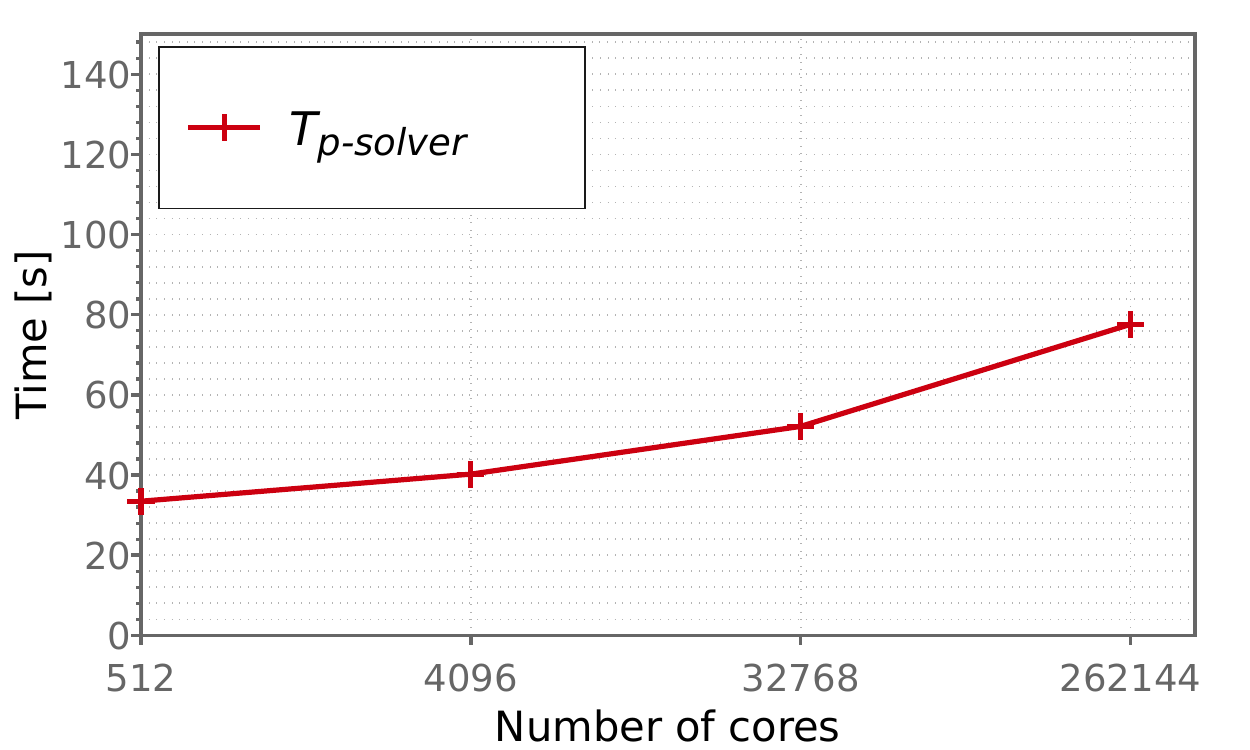}}
 	    \caption{Time measurements}
		{\includegraphics[width=0.4999\textwidth]{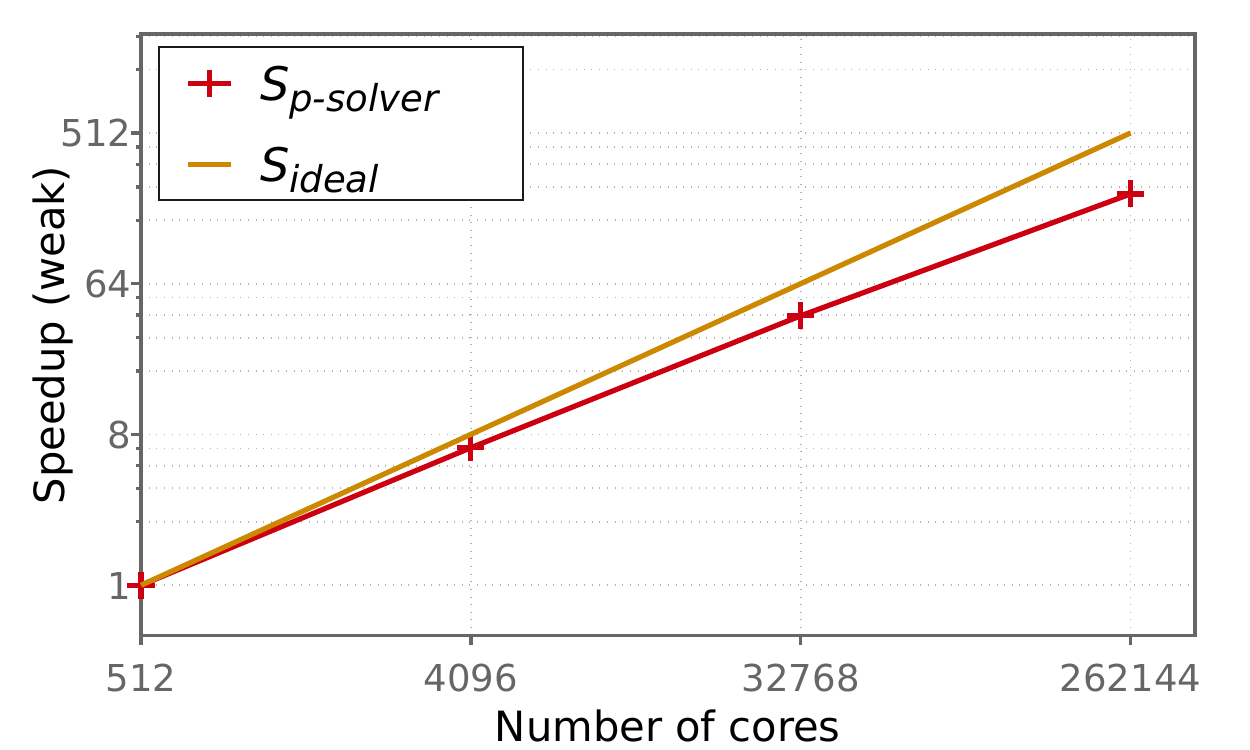}}	
		\caption{Speedup relative to 512 cores}
		\label{fig:TimeAndSpeedup}
	\end{subfigure}	
	\vskip\baselineskip
	\begin{subfigure}[t]{\textwidth}
		\centering
		\begin{tabular}{rrrccc}
			\centering
			\begin{tabular}{llllccc}
				\toprule
				Procs  & Refs  & Num.Elems  & DoFs & Newton Its. &Total Lin.Its.&Lin.Its. \cref{eq:LagrangianNewtonSystem1}\\ 
				\midrule
				512 & 2 & \num{77824} & {44730} & 68& \num{2080} & 394  \\
				\num{4096} & 3 & \num{622592} & \num{334158} & 68& \num{2458} & 472  \\
				\num{32768} & 4 & \num{4980736} & \num{2581014}& 68& \num{2606} & 509  \\
				\num{262144}& 5 & \num{39845888} & \num{20283942} & 68& \num{2912} & 577  \\
				\bottomrule
			\end{tabular}
		\end{tabular}
		\caption{Iteration counts for one optimization step of the solver used to obtain $u_p$, lines 11-15 in~\cref{alg:optimizationAlgorithm}.}
		\label{tab:iter_counts}
	\end{subfigure}
	\caption{Weak Scaling: Results for the first optimization step. Accumulated wallclock time for all $p$-levels and speedup relative to 512 cores are shown. The number of Newton steps across several levels of refinement, as well as the linear solver iterations are presented in relation to the number of tetrahedrons per refinement level and the corresponding DoFs in the discretization of \cref{eq:DiscreteSaddlePointProblem}.}
	\label{fig:weak_scaling}
\end{figure}

A 3d computational grid with 2 levels of refinement is used as an initial measurement in order to optimize the number of cores used at the finest level. 
The wallclock times, speedup, and iteration counts are shown in  \cref{fig:weak_scaling}.
An eight-fold increase in the number of cores is performed for each level of refinement
Results are presented for the solution of the nonlinear system of equations given in \cref{eq:optimality_system} via its linearization in \cref{eq:LagrangianNewtonSystem1,eq:LagrangianNewtonSystem2}. 
This system is solved using Newton's method with a BiCGStab as a solver for the underlying linearization. 
The linear solver is set to absolute and relative error reductions of \num{1e-10}  and  \num{1e-16}, respectively. 
It is preconditioned by a geometric multigrid method with 3 pre- and post smoothing steps via a Gauss-Seidel smoother within a V-cycle. 
An LU factorization solves the base level gathered in a single core.

We measure the accumulated times and iteration counts for the routines in lines 11-15 of~\cref{alg:optimizationAlgorithm} for one optimization step. 
This can be understood as the time it takes to assemble the linearization, initialize the grid hierarchy necessary for the geometric multigrid preconditioner, and apply the linear solver until convergence within each call to the Newton's solver.
This is done for each value of $p$ starting at $p_\mathrm{init}$ up to $p_\mathrm{max}$ with $p_\mathrm{inc}$ intervals as explained in~\cref{sec:algorithm}. 
The time measurement starts for every optimization at $u_{p_\mathrm{init}}$ and ends once the corresponding Newton solver for $u_{p_\mathrm{max}}$ has converged.
The speedup (b) is presented relative to the base measurement with 512 cores, and in (c) the iteration counts are shown in relation to the number of DoFs and tetrahedral elements.
The column of the total linear iterations includes all the necessary calls to the linear solver used within the linearization. 
As shown in \cref{eq:SchurComplementReducedSystem}, for each solution of the linear system of equations it is necessary to solve $m+2$ times with $A^{-1}$. 
These include one time for the rhs in the second equation of \cref{eq:SchurComplementReducedSystem}, and $m$ for the computation of $S$. 
Addtionally, the first equation of \cref{eq:SchurComplementReducedSystem} has to be solved for $\delta_u$, whose iteration counts are shown, individually, in the rightmost column of \cref{fig:weak_scaling}(c).

It can be seen that good scalability results are obtained for up to \num{262144} cores.
The communication costs impose a time overhead significantly lower to the very large increase of the number of DoFs.
Altogether, the results show the need for using numerical solvers with grid independent convergence.
Recall that our target is to use the solution of the $p$-Laplace relaxed problem for the highest value of $p$, i.e. $p_\mathrm{max}$, as a deformation field to generate a series of shape iterates. 
Moreover, we do this by solving the same problem for lower values of $p$, in order to have a good initial guess as we approach the maximum $p$. 
The latter fact is necessary, since with each increment of $p_\mathrm{inc}$, our problem becomes more nonlinear, implying it becomes more difficult to solve, particularly without a good initial guess.
For the given settings, $p_\mathrm{init}=2.0$ to $p_\mathrm{max}=4.1$ and an increment of $p_\mathrm{inc}=0.19$, \cref{alg:NewtonSolver} has to be called thirteen times.
Newton's method has to call the linear solver for each of these $p$ values.
Therefore, there is an evident need for an efficient, fast, and computationally cheap preconditioner which allows for grid-size independent bounds on the convergence rate of the iterative methods.
This is possible with the geometric multigrid method.
One of the downsides is that this preconditioner requires a base level computational mesh that describes a  geometry that can be represented by a grid hierarchy, see \cite{Reiter2014gmg}, which implies that care must be taken during the generation of the grid. 
Nevertheless, it is a very effective approach towards solving for $u_p$ with increments of the $p$ value. 
The results in \cref{fig:weak_scaling} show that the $p$-relaxed problem becomes inexpensively solvable. 
Additionally, the benefits of the multigrid preconditioner are evident by noticing how the Newton's method is perfectly scalable in the number of steps needed for all refinement levels, as well as in the slight increase in linear solver iterations between the initial and final runs.
As seen in the table, even when the number of DoFs increases by three orders of magnitude, the timings and iteration counts are bound by the preconditioner.

\begin{figure}[!htbp]
	\centering
	\includegraphics[width=0.8\textwidth]{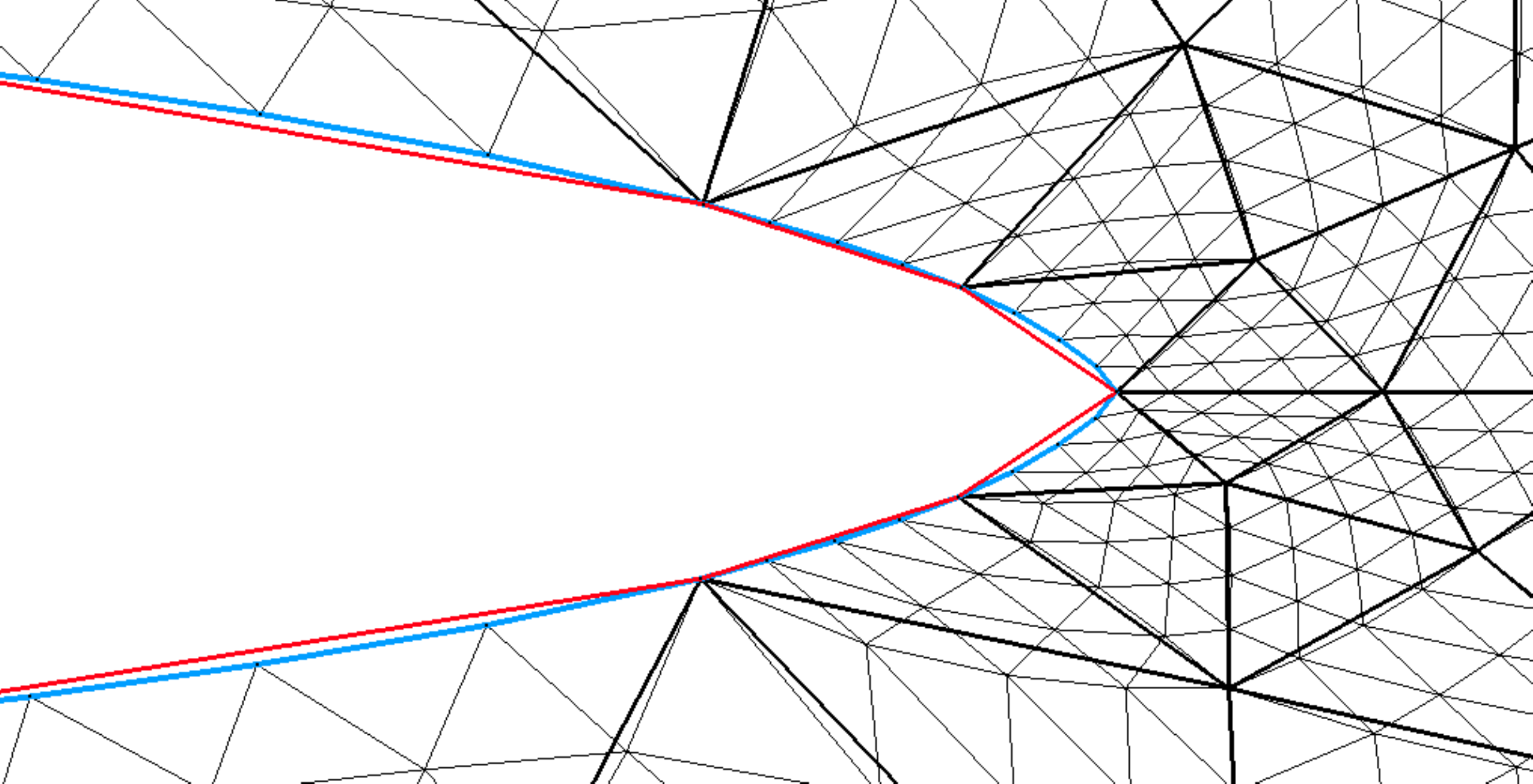}
	\caption{The base level and the second level of refinement are compared for a 2d simulation for the last step before convergence. The macro and refined triangular elements are shown in bold and thin black lines, respectively. Given that the deformation field is restricted and applied throughout the grid hierarchy, the coarsest grid is an interpolation of the finest.}
	\label{fig:coarse-grid}
\end{figure}

In order to preserve numerical scalability across all optimization steps, it is necessary to apply the deformation field across the complete grid hierarchy. 
This is shown in \cref{fig:coarse-grid}, where the base level is compared to the finest grid with two refinements. 
It is visible how $u_p$ is restricted and applied to all levels, therefore generating an optimal coarse grid. Given that this implies an interpolation of the vector field, and that by definition the obstacle's surface on the coarsest grid has less nodes than the upper levels, there is a slight mismatch between the two grids. 
However, this has no detrimental effects nor adds more computational complexity to the shape optimization scheme. Our scheme works on arbitrary Lipschitz shapes. 
Therefore, it is not necessary to incorporate extra geometric information to the grid hierarchy.

Overall, good weak scalability results were obtained for up to 39 million elements. 
This represents an increase of three orders of magnitude, both in tetrahedrons as in DoFs, with a slight increase in the necessary computational work in terms of linear iterations.
Although the performance dropped marginally from the ideal case, the wallclock times and speedup show that the numerical scheme we propose for the solution of the $p$-Laplace relaxed problem could be used for problems with large numbers of DoFs, corresponding to real-world industrial applications.

\section{Conclusion}
\label{sec:conclusion}
In this work we presented a steepest descent shape method based on $W^{1,p}$ approximations of $W^{1,\infty}$ for shape optimization problems with PDE and fixed-dimensional geometric constraints.
We demonstrated that the algorithms works for general Lipschitz shapes since deformations allow singularities in the surface to be smoothed or newly generated.
Furthermore, we incorporated fixed-dimensional constraints together with the PDE constraints into the optimization algorithm via a Schur Complement approach.
Compared to approximate algorithms such as the penalty based and augmented Lagrangian approaches, we demonstrated a significant gain of robustness in the treatment of geometric constraints over the optimization steps.
Additionally, this work addressed line search schemes for the steepest descent direction in $W^{1,\infty}$, which -in  contrast to the aforementioned Hilbert space methods- is here a nonlinear problem. Moreover, this problem lives on the solution manifold of the nonlinear geometric constraints posing a non-convex set in general.

The essential part of this work was to investigate the application of geometrical multigrid preconditioners on hierarchical grid structures without needing any further information, such as curvature based on spline surfaces.
It was demonstrated that, via the shape optimization, a body-fitted hierarchical grid structure is found for the optimal shape.
Our numerical studies indicated that under this circumstances the multigrid preconditioner features a mesh-independent solver for the deformation subproblem.
As a consequence, we were able to demonstrate that the proposed method exhibits weak scalability up to \num{262144} CPU cores of the distributed-memory system Hawk at HLRS.

\subsection*{Acknowledgements}
{	
	\noindent
	The current work is part of the research training group ``Simulation-Based Design Optimization of Dynamic Systems Under Uncertainties'' (SENSUS) funded by the state of Hamburg under the aegis of the Landesforschungsförderungs-Project LFF-GK11.
	
	\noindent
	The authors acknowledge the support by the Deutsche Forschungsgemeinschaft (DFG) within the Research Training Group GRK 2583 ``Modeling, Simulation and Optimization of Fluid Dynamic Applications''.

	\noindent
	Computing time on the national supercomputer HPE Apollo Hawk at the High Performance Computing Center Stuttgart (HLRS) under the grant ShapeOptCompMat (ACID 44171, Shape Optimization for 3d Composite Material Models) is gratefully acknowledged.
}
\printbibliography

@book{sokolowski1992,
    author = {Sokolowski, J. and Zol{\'e}sio, J.  -P.},
    title = {Introduction to Shape Optimization},
    subtitle = {Shape Sensitivity Analysis},
    year = {1992},
    publisher = {Springer, Berlin, Heidelberg},
    DOI = {10.1007/978-3-642-58106-9}
}

@article{pironneau1973optimum,
    title={On optimum profiles in Stokes flow},
    author={Pironneau, O.},
    journal={Journal of Fluid Mechanics},
    volume={59},
    number={1},
    pages={117--128},
    year={1973},
    publisher={Cambridge University Press},
    DOI={10.1017/S002211207300145X}
}

@book{zolesio2011,
    author = {Delfour, M. C. and Zol{\'e}sio, J. -P.},
    title = {Shapes and geometries: metrics, analysis, differential calculus, and optimization},
    year = {2011},
    publisher = {SIAM},
    DOI = {10.1137/1.9780898719826}
}

@article{allaire2004structural,
    title={Structural optimization using sensitivity analysis and a level-set method},
    author={Allaire, G. and Jouve, F. and Toader, A.-M.},
    journal={Journal of computational physics},
    volume={194},
    number={1},
    pages={363 -- 393},
    year={2004},
    publisher={Elsevier},
    DOI={10.1016/j.jcp.2003.09.032}
}

@inproceedings{allaire2020,
  title={Chapter 1 - Shape and topology optimization},
  author={Allaire, G. and Dapogny, C. and Jouve, F.},
  Booktitle = {Geometric Partial Differential Equations - Part II},
  year = {2021},
  Series = {Handbook of Numerical Analysis},
  volume = {22},
  publisher = {Elsevier},
  pages = {1 -- 132},
  DOI = {10.1016/bs.hna.2020.10.004}
}

@article{schulz2016,
    title={Computational Comparison of Surface Metrics for PDE Constrained Shape Optimization},
    author={Schulz, V. and Siebenborn M.},
    journal={ Computational Methods in Applied Mathematics},
    volume={16},
    number={3},
    pages={485--496},
    year={2016},
    publisher={De Gruyter},
    DOI={10.1515/cmam-2016-0009}
}

@article{onyshkevych2020,
    title={Mesh Quality Preserving Shape Optimization Using Nonlinear Extension Operators},
    author={Onyshkevych, S. and Siebenborn, M.},
    journal={Journal of Optimization Theory and Applications},
    volume={189},
    pages={291--316},
    year={2020},
    DOI={10.1007/s10957-021-01837-8}
}

@article{haubner2021,
    title={A Continuous Perspective on Shape Optimization via Domain Transformations},
    author={Haubner, J. and Siebenborn, M. and Ulbrich, M.},
    journal={SIAM},
    volume={43},
    number={3},
    pages={A1997--A2018},
    year={2021},
    DOI={10.1137/20M1332050}
}

@article{deckelnick2021,
    AUTHOR = {Deckelnick, K. and Herbert, P. and Hinze, M.},
	Title = {A novel $W^{1,\infty}$ approach to shape optimisation with Lipschitz domains.},
	Year = {2021},
	archivePrefix={arXiv},
	eprint={2103.13857}
}

@article{mueller2021,
    title={A novel p-harmonic descent approach applied to fluid dynamic shape optimization},
    author={M\"uller, P. M. and K\"uhl, N. and Siebenborn, M. and Deckelnick, K. and Hinze, M. and Rung, T.},
    journal={Struct Multidisc Optim}, 
    year={2021},
    DOI={10.1007/s00158-021-03030-x}
}

@inproceedings{brandenberg2011,
    title={Advanced Numerical Methods for PDE Constraint Optimization with Application to Optimal Design in Navier Stokes Flow},
    author={Brandenburg, C. and Lindemann, F. and Ulbrich, M. and Ulbrich, S.},
    booktitle={Constrained Optimization and Optimal Control for Partial Differential Equations},
    editors={Kunisch, K. and Sprekels, J. and Leugering, G. and Tr\"olzsch, F.},
    year={2011},
    series={International Series of Numerical Mathematics},
    volume={160},
    pages={257--275},
    publisher={Brinkh\"aser, Basel},
    DOI={10.1007/978-3-0348-0133-1_14}
}

@book{mohammadi2010,
    title={Applied Shape Optimization for Fluids},
    author={Mohammadi, B. and Pironneau, O.},
    publisher={Oxford University Press},
    year={2009},
    DOI={10.1093/acprof:oso/9780199546909.001.0001}
}

@misc{schiela2019composite,
          author = {Ortiz, J. and Schiela, A.},
           title = {A composite step method for equality constrained optimization on manifolds},
            year = {2019},
             url = {https://eref.uni-bayreuth.de/47908/},
}

@article{schiela2021sqp,
  doi = {10.1137/20m1341325},
  year = {2021},
  publisher = {Society for Industrial {\&} Applied Mathematics ({SIAM})},
  volume = {31},
  number = {3},
  pages = {2255--2284},
  author = {Schiela, A. and Ortiz, J.},
  title = {An {SQP} Method for Equality Constrained Optimization on Hilbert Manifolds}
}

@article{mohammadi2004,
    author={Mohammadi, B. and Pironneau, O.},
    title={Shape Optimization in Fluid Mechanics},
    journal={Annual Review of Fluid Mechanics},
    volume={36},
    pages={255--279},
    year={2004},
    doi={10.1146/annurev.fluid.36.050802.121926}
}

@article{ulbrich2003,
    author={Ulbrich, M.},
    title={Constrained Optimal Control of Navier-Stokes Flow by Semismooth Newton Methods},
    journal={Systems \& Control Letters},
    volume={48},
    pates={297--311},
    year={2003},
    doi={10.1016/S0167-6911(02)00274-8}
}

@article{hinze2001,
    author={Hinze, M. and Kunisch, K.},
    title={Second order methods for optimal control of time-dependent fluid flow},
    journal={SIAM},
    volume={40},
    pages={925-946},
    year={2001},
    doi={/10.1137/S0363012999361810}
}

@article{ishii2005limits,
    title={Limits of solutions of p-Laplace equations as p goes to infinity and related variational problems},
    author={Ishii, H. and Loreti, P.},
    journal={SIAM journal on mathematical analysis},
    volume={37},
    number={2},
    pages={411 -- 437},
    year={2005},
    publisher={SIAM},
    DOI={10.1137/S0036141004432827}
}

@misc{pinzon2021,
    author={Pinzon, J. and Siebenborn, M.},
    title={Fluid dynamic shape optimization using self-adapting nonlinear extension operators with multigrid preconditioners},
    year={2021},
    archivePrefix={arXiv},
    eprint={2108.07788v1}
}

@article{Reiter2014gmg,
	Author = {Reiter, S. and Vogel, A. and Heppner, I. and Rupp, M. and Wittum, G.},
	Journal = {Comp. Vis. Sci.},
	Number = {4},
	Pages = {151-164},
	Publisher = {Springer},
	Title = {A massively parallel geometric multigrid solver on hierarchically distributed grids},
	Volume = {16},
	Year = {2013}
}

@article{Vogel2014ug4,
	Author = {Vogel, A. and Reiter, S. and Rupp, M. and N\"agel, A. and Wittum, G.},
	Journal = {Comp. Vis. Sci.},
	Number = {4},
	Pages = {165-179},
	Title = {{UG} 4: A novel flexible software system for simulating {PDE} based models on high performance computers},
	Volume = {16},
	Year = {2013}
}

@ARTICLE{GMSH,
	AUTHOR = {Geuzaine, C. and Remacle, J.-F.},
	TITLE = {Gmsh: {A} 3-{D} finite element mesh generator with built-in
	pre- and post-processing facilities},
	JOURNAL = {International Journal for Numerical Methods in Engineering},
	VOLUME = {79},
	YEAR = {2009},
	NUMBER = {11},
	PAGES = {1309--1331},
	ISSN = {0029-5981},
	DOI = {10.1002/nme.2579},
}

@misc{PARMETIS,
	title={Parmetis, Parallel graph partitioning and sparse matrix ordering library},
	author={Karypis, G. and Schloegel, K. and Kumar, V.},
	url={http://glaros.dtc.umn.edu/gkhome/metis/parmetis/overview},
	urldate={2020-03-06},
	version={4.03},
	year={2013}
}

@book{Hackbusch85,
	Author = {Hackbusch, W.},
	Publisher = {Springer},
	Title = {Multi-grid methods and applications},
	Volume = {4},
	Year = {1985}}

@webpage{VTK,
	Url = {vtk.org}
}

@inproceedings{geiser2021aggregated,
	title={Aggregated Formulation of Geometric Constraints for Node-based Shape Optimization with Vertex Morphing (Eurogen 2021)},
	author={Geiser, A. and Antonau, I. and Bletzinger, K.-U.},
	booktitle={14th ECCOMAS Thematic Conference on Evolutionary and Deterministic Methods for Design, Optimization and Control},
	editors={Gauger, N. and Giannakoglou, K. and Papadrakakis, M. and Periaux, J.},
	year={2021},
	volume={14},
	pates={80--94},
	publisher={Eccomas Proceedia},
	DOI={10.7712/140121.7952.18383}

}

@article{uzawa1958iterative,
	title={Iterative methods for concave programming},
	author={Uzawa, Hirofumi},
	journal={Studies in linear and nonlinear programming},
	volume={6},
	pages={154--165},
	year={1958}
}

@article{arora1991,
	title={Multiplier methods for engineering optimization},
	author={Arora, J.S. and Chahande, A.I. and Paeng J.K.},
	journal={International Journal for Numerical Methods in Engineering},
	volume={32},
	number={7},
	pages={1485--1525},
	year={1991},
	publisher={Wiley},
	DOI={10.1002/nme.1620320706}
}

@misc{gmgshapeopt2021,
	title={PLaplaceOptim},
	author={Pinzon, J. and Siebenborn, M.},
	howpublished ={\url{http://www.github.com/multigridshapeopt}},
	urldate={2021-08-17},
	version={1.0},
	year={2021}
}

\end{document}